\documentclass[preprint,review,10pt,authoryear]{elsarticle}

\usepackage{enumerate,mathrsfs,setspace,geometry,epstopdf,bm,arydshln,amsmath,multirow,amsthm,color}
\usepackage{longtable,ulem,extpfeil,extarrows,graphicx,cases,subfigure,booktabs}%
\usepackage{amssymb,lscape,diagbox}
\usepackage{graphicx}
\usepackage{verbatim}
\usepackage[ruled,vlined,linesnumbered]{algorithm2e}
\usepackage[colorlinks,citecolor=black]{hyperref}
\usepackage{natbib}
\usepackage{amsmath}
\usepackage{tabularx}
\usepackage{multirow}
\usepackage{multicol}
\usepackage{longtable}
\usepackage{ltxtable}
\usepackage{appendix}
\usepackage{color}
\usepackage{subfigure}
\usepackage{rotating}
\usepackage{geometry}
\usepackage{enumerate,diagbox}
\usepackage{bm}
\usepackage{tabu}
\newtheorem{lemma}{Lemma}[section]
\newtheorem{proposition}{Proposition}[section]
\newtheorem{theorem}{Theorem}[section]
\newtheorem{definition}{Definition}[section]
\newtheorem{remark}{Remark}[section]
\DeclareMathOperator*{\esssup}{ess\,sup}
\DeclareMathOperator*{\sta}{s.t.}

\allowdisplaybreaks[2]

\newcommand{\ignore}[1]{}

\begin{document}

\begin{frontmatter}


\title{Multi-period facility location and capacity planning under $\infty$-Wasserstein joint chance constraints in humanitarian logistics}

\author[a1]{Zhuolin Wang}\ead{wangzl17@mails.tsinghua.edu.cn}
\author[a1]{Keyou You}\ead{youky@tsinghua.edu.cn}
\author[a3]{Zhengli Wang}\ead{zhlwang@bjtu.edu.cn}
\author[a3]{Kanglin Liu\corref{cor1}}\ead{klliu@bjtu.edu.cn}

\cortext[cor1]{Corresponding author}

\address[a1]{Department of Automation and BNRist, Tsinghua University, Beijing 100084, P. R. China}
%\address[a2]{Institute of Transportation System Science and Engineering, Beijing Jiaotong University, Beijing 100044, China}
\address[a3]{School of Traffic and Transportation, Beijing Jiaotong University, Beijing 100044, P. R. China}

\begin{abstract}
% Determine the capacity of facilities in humanitarian logistics.
The key of the post-disaster humanitarian logistics (PD-HL) is to build a good facility location and capacity planning (FLCP) model for delivering relief supplies to affected areas in time. To fully exploit the historical PD data, this paper adopts the data-driven distributionally robust (DR) approach and proposes a novel multi-period FLCP  model under the $\infty$-Wasserstein joint chance constraints (MFLCP-W). Specifically, we sequentially decide locations  from a candidate set to build facilities with supply capacities,  which are expanded if more economical, and  use a \textit{finite} number of historical demand samples in chance constraints to ensure a high probability of on-time delivery. To solve the MFLCP-W model, we equivalently reformulate it as a mixed integer second-order cone program and then solve it by designing an effective outer approximation algorithm with two tailored valid cuts. Finally, a case study under hurricane threats shows that  MFLCP-W outperforms its counterparts  in the terms of the cost and service quality, and that our algorithm converges significantly faster than the commercial solver CPLEX 12.8 with a better optimality gap.   %other data-driven single-period FCLP and MFLCP models.

%minimize the worst-case cost and satisfy the worst-case demand supply constraints over a $\infty$-Wasserstein ball, which is centered at the observed empirical distribution, with the ball radius reflecting its confidence level. Different from existing DR facility location and capacity planning models, our multi-period model accounts for the changeable disaster situations and the human sufferings. Moreover, we propose a novel penalty cost function to quantify human sufferings. In sharp contrast to most intractable DR joint chance constrained problems, the proposed model over the $\infty$-Wasserstein ball

\end{abstract}

\end{frontmatter}

\section{Introduction}
\label{Section:intro}
Disasters have greatly threatened to human societies.  To reduce their impact, it is essential to deliver relief supplies to affect areas in time after the outbreak of a disaster, which is the main task of the post-disaster humanitarian logistics (PD-HL). This usually requires to build a good facility location and capacity planning (FLCP)  model and provides an allocation strategy for transporting relief supplies.
Such a problem has attracted an increasing attention, and involves at least the following major challenges.

The first is how to exploit changing features of affected areas in an FLCP model \citep{Jose2013On}. In the literature, they usually use initial information for prediction and then make a static plan of a finite-time horizon, see e.g. \citet{rawls2010pre,li2010continuum,chen2016network,elcci2018chance,O2018Chance,ni2018location,liu2019distributionally,velasquez2020prepositioning}.  Clearly, the resulting solution cannot adapt to the PD data that is sequentially collected from affected areas. With more relief supplies transported to the affected areas, there is no reason to keep our initially-made conservative plan. To resolve it, we propose a multi-period FLCP (MFLCP) model where the decisions in later periods will be modified based on the results of former periods. %only if the decisions in former periods have been revealed. 
Specifically, both facility locations and their capacities are to be jointly adjusted by adapting to spatial-temporal evolution of the affected areas. Though the models in \citet{blood2014TRE,Charles2016Designing,ALEM2020} also involve a dynamical PD-HL setting, they mainly focus on the location of facilities without adjusting their capacities. In fact, the capacity design can reduce both the capacity redundancy and the total operational cost  \citep{melkote2001capacitated,filippi2021single,saif2021data}.

The second is how to jointly minimize costs from both the supply-side and demand-side. The existing FLCP models essentially only minimize the supply-side cost, including the set-up cost of new facilities, the capacity expansion cost, the transportation cost and other logistic costs. Differently, this work further considers the demand-side penalty cost of human suffering from the lack of critical supplies over a period \citep{Jose2013On}. Quantifying human sufferings is of independent interest and has become a promising research direction. Two main ideas use the deprivation cost function (DCF) \citep{Jose2013On} and the deprivation level function (DLF) \citep{Wang2017poms}, respectively.  A DCF evaluates the  willingness-to-pay for relief supplies from the perspective of welfare economics  \citep{Jose2013On,holguin2016econometric} and admits an economical value that can be easily added to  the supply-side cost. However, it is greatly affected by local economic levels and individual incomes, rendering it difficult to compare among different cases \citep{shao2020research}. Moreover, the DCF only collects data for a short planning horizon right after a disaster (say 48 hours \citep{holguin2016econometric}), and thus is not suitable for the practical relief process that may last for weeks.  In a comparison, a DLF adopts a numerical rating scale method to quantify the degree of human sufferings \citep{Wang2017poms,Wang2019Augmenting}. Yet, its degree value is dimensionless and it is not easy to compare/add/subtract with the supply-side cost. Inspired by these observations, we evaluate human sufferings by proposing a demand-side penalty cost using the weighted unmet demand where the weighting factor is nonlinear and monotonically increasing with respect to the deprivation time. Such a simple idea can not only overcome the shortcomings of DCF and DLF, but also inherit their advantages, e.g., it can be easily used in the long-time horizon case. Note that the penalty-cost based methods  with static weights have been extensively adopt in HL  \citep{rawls2010pre,kelle2014decision,khayal2015model,shao2020research}.

The last is how to satisfy the  demand of relief supplies in the uncertain environment. While the demand is unpredictable and subject to large uncertainty, the classical robust model aims to satisfy all demands over a given set \citep{qiaofengli}, leading to a very conservative solution \citep{liu2019distributionally}. The stochastic model assumes that the demand follows a probability distribution and in practice is solved by the sample average approximation (SAA) approach  \citep{O2018Chance,elcci2018chance}.  To take advantage of both models,  the distributionally robust (DR) approach has been adopt in \citet{liu2019distributionally,saif2021data,jiang2021distributionally} under the assumption that the distribution of demand belongs to a set of probability distributions, over which an optimal solution is to be found in some worst-case sense, e.g., a solution maximizes the worst-case probability of the on-time delivery.  

% In fact, there has been a surge interest in the application of the DR problems, such as the knapsack problem \citep{cheng2014distributionally}, the vehicle routing problem   \citep{ghosal2020distributionally} and supply chain network problem \citep{karimi2021biobjective}. 

A key issue of the DR FLCP model is the construction of a good ambiguity set, which should be large enough to contain the true distribution with a high probability, but cannot be too ``large" to avoid overly-conservative decisions.
To our knowledge, the DR FLCP model mainly uses a moment-based ambiguity set of probability distributions with some specified moment constraints \citep{liu2019distributionally}. Under \textit{exact} first- and second-order moments, the DR FLCP model can be reformulated as a tractable conic optimization problem \citep{khodaparasti2018multi}. Noting that the moment mismatch is unavoidable, \citet{yang2021multi} further incorporate moment uncertainty in their FLCP model and obtain a strongly NP-hard problem.  In this work,  the ambiguity set is a data-driven $\infty$-Wasserstein ball \citep{kantorovich1958space}  centered at the empirical distribution of samples of demands with an appropriately selected radius, which reflects our confidence level of samples, e.g., the smaller the radius, the larger the confidence, and captures the slowly time-varying nature of the demand distribution. 

To sum up, we propose a novel MFLCP model with $\infty$-Wasserstein joint chance constraints (MFLCP-W) and a new demand-side penalty cost function to quantify human sufferings. As the DR $1$-Wasserstein chance constrained models in \citet{xie2021distributionally, chen2018data},  MFLCP-W is an infinite-dimensional problem and is typically difficult to solve. In contrast to their intractable bilinear reformulation, we show that  the MFLCP-W model can be equivalently solved by a finite mixed integer second-order conic program (MISOCP) and then solve it by designing an outer approximation (OA) algorithm \citep{1987An}  with two valid cuts, which  significantly outperforms the commercial  CPLEX 12.8 in terms of the computational time and optimality gap.  Moreover, a case study under hurricane threats shows that  MFLCP-W outperforms its counterparts  in the terms of the cost and service quality.

%, which can obtain the optimal solution within finite number of iterations \citep{BONAMI2008186}. 

%To the best of our knowledge, we are the first to adopt the joint chance constraints over the $\infty$-Wasserstein set in the MFLCP model. The most relevant paper that involves the Wasserstein set in a facility location problem is \cite{saif2021data}. They propose a single-period model and require solution to remain feasible for all demand realizations, which is more conservative than ours with joint chance constraints. 
%这里那些引用我删掉了，感觉好像没啥必要，强调一下我们这个算法的有效性就行了吧？
%OA algorithm is proposed by \cite{1987An} to solve mixed integer nonlinear programs (MINLP), and is able to obtain the optimal solution within finite number of iterations if the continuous relaxation of MINLP is proved to be convex \citep{BONAMI2008186}. Due to its computational tractability and theoretical guarantee, the OA algorithm has been widely used in applications such as closed-loop supply chain \citep{LIU2018244,qiaofengli}, supply chain network design \citep{vanaEJORsubmodular201846,MAI2020874} and production management \citep{TORKAMAN2020105019}.
%Thus, the reformulated MISOCP, which is polynomial time solvable, outperforms  other DR chance constrained models with left-hand-side uncertainty which are proved to be NP-hard \citep{hanasusanto2017ambiguous}.

The reminder of this paper is organized as follows. In Section \ref{section:model}, we present a detailed problem statement and model formulation of the MFLCP-W model. Section \ref{section:reformulation} equivalently reformulates the proposed model into an MISOCP. Section \ref{section:oa} introduces a tailored OA algorithm to solve the reformulated model. Numerical results for a case study are conducted in Section \ref{section:numerical}. Finally, we conclude the paper and propose several future directions in Section \ref{section:conclusion}.

\section{Problem formulation and model development} \label{section:model}
%In this section, we describe our MFLCP-W model in details.  

\begin{figure}[htbp]
 \centering
 \includegraphics[width=0.8\linewidth]{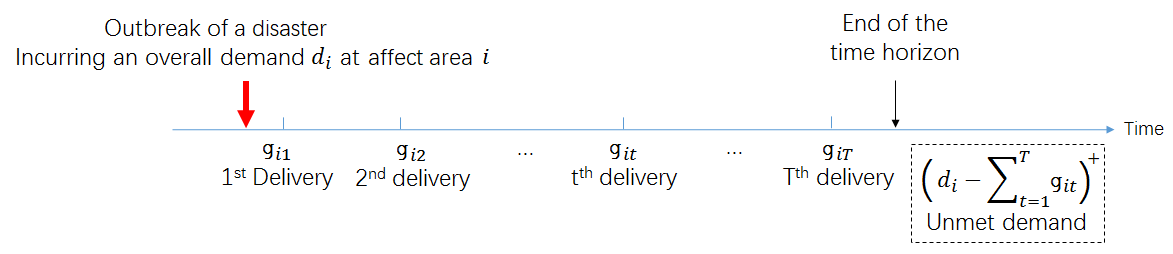}
 \caption{Sketch of fulfilling the demand  process.}
 \label{fig:sketch}
\end{figure}

Consider a relief network that is composed of  a set of demand sites (affected areas) and a set of candidate locations for constructing facilities, which are denoted as $[I]$ and $[J]$ respectively, and  each demand site $i\in[I]$ incurs an unknown demand $d_i$ of relief supplies.  The objective of PD-HL is to decide when and how to select locations from $[J]$ to construct facility with their supply capacities, and design an allocation strategy for transporting relief supplies from facilities to  demand sites in an ``optimal" way. A key ingredient of this work is that the random demand $d_i$ is  not directly observed and can only be estimated from a sequence of historical demand samples which poses significant challenge in the PD-HL.  Moreover, we focus on the multi-period planning for a finite-time horizon $T$ and allow unmet demands at the end of the horizon. Specifically, let $g_{it}$ be the total of relief supplies delivered to demand site $i$ at period $t$. The unmet demand of this site is expressed as $\left(d_i-\sum_{t \in [T]}g_{it}\right)^{+}$ where $\left(\cdot\right)^{+}$ returns the non-negative part of its argument.  See Figure \ref{fig:sketch} for a simple illustration on how to (partially) meet the demand $d_i$ in a finite-time horizon.

Without loss of generality, we adopt the assumption in \cite{Shulman1991An} and \cite{yu2020multistage} that the constructed facilities will not get destroyed or closed, and  their capacities of each facility will not be reduced before the planning horizon $T$. By convention, a dummy period $0$ is introduced for the sake of completeness and the capacities of all facilities are initially set to zero.
%only constructed facilities will not get destroyed or closed during the planning horizon. (a) Once built, the facility will not be destroyed. (b) $S_{j0}=0$. Then the standard multi-period location and capacity planning (MLCP) model aims to plan distribution of the facilities and assign proper capacity to them to satisfy the medical demand of all affected areas with the objective of minimizing the total logistic cost, i.e.,
Notations of this paper are summarized as follows. 

% \noindent {\bf Indices}

% \vspace{0.3cm}\begin{tabular}{ll}
% %\multicolumn{2}{l}{Sets and indices}\\
%  $j$ & The location index of the candidate facility,  e.g.,  $j\in [J]=\{1,\ldots,J\}$.\\
%  $i$ & The index of the demand site,  e.g.,  $i\in [I]=\{1,\ldots,I\}$.\\
%  $t$& The index of planning time period,  e.g., $t\in [T]=\{1,\ldots,T\}$.
%  \end{tabular}
 
\renewcommand{\arraystretch}{1} %¿ØÖÆ±í¸ñÐÐ¸ßµÄËõ·Å±ÈÀý
\begin{table}[!h]
	\caption{Indices}
	%	\fontsize{8}{10}\selectfont
	\begin{tabu} to 1\textwidth {X[1,l] X[6,l] }
		\toprule[1 pt]
		Notation & Description  \\ \midrule
		$j$ & The location index of the candidate facility,  e.g.,  $j\in [J]=\{1,\ldots,J\}$.\\
        $i$ & The index of the demand site,  e.g.,  $i\in [I]=\{1,\ldots,I\}$.\\
         $t$& The index of planning time period,  e.g., $t\in [T]=\{1,\ldots,T\}$.\\
		\bottomrule[1 pt]		
	\end{tabu}
	\label{indices}	
\end{table}

\renewcommand{\arraystretch}{1} %¿ØÖÆ±í¸ñÐÐ¸ßµÄËõ·Å±ÈÀý
\begin{table}[!h]
	\caption{Parameters given as {\it prior}}
	%	\fontsize{8}{10}\selectfont
	\begin{tabu} to 1\textwidth {X[1,l] X[6,l] }
		\toprule[1 pt]
		Notation & Description  \\ \midrule
		$f_{jt}$& The facility set-up cost in location $j$ in period $t$.\\
        $a_{jt}$& The unit cost of expanding capacity in facility  $j$ in period $t$. \\
        $q_t$& The  limit of supply capacity of every facility in period $t$. \\
        $c_{ji}$& The unit  cost of transporting relief supplies from location $j$ to demand site $i$.\\
		\bottomrule[1 pt]		
	\end{tabu}
	\label{param_prior}	
\end{table}
 
\renewcommand{\arraystretch}{1} %¿ØÖÆ±í¸ñÐÐ¸ßµÄËõ·Å±ÈÀý
\begin{table}[!h]
	\caption{Random variables and their samples}
	%	\fontsize{8}{10}\selectfont
	\begin{tabu} to 1\textwidth {X[1,l] X[6,l] }
		\toprule[1 pt]
		Notation & Description  \\ \midrule
		$d_i$ & The random demand in the demand site $i \in [I]$. \\
        $\hat{d}_i^h$ & An estimate/sample of $d_i$ that was obtained from the $h$-th historical observation. \\
		\bottomrule[1 pt]		
	\end{tabu}
	\label{random_var_sample}	
\end{table}

\renewcommand{\arraystretch}{1.2} %¿ØÖÆ±í¸ñÐÐ¸ßµÄËõ·Å±ÈÀý
\begin{table}[!h]
	\caption{Decision variables}
	%	\fontsize{8}{10}\selectfont
	\begin{tabu} to 1\textwidth {X[1,l] X[10,l] }
		\toprule[1 pt]
		Notation & Description  \\ \midrule
		$X_{jt}$ &Binary variable from $\{0, 1\}$  and $X_{jt}=1$ if a facility is constructed on location $j$  in period $t$.\\
        $S_{jt}$ & The discretized supply capacity of facility $j$ in period $t$, e.g., $S_{j0}=0$. \\
        $Y_{jit}$& Continuous variable in $[0, 1]$ to denote the percentage of $d_i$ served by facility $j$ in period $t$.\\
        $g_{it}$ &The total of relief supplies delivered to demand site $i$ in period $t$, i.e. $g_{it}=\sum_{j\in [J]}Y_{jit}d_i$.\\
		\bottomrule[1 pt]		
	\end{tabu}
	\label{dec_var}	
\end{table}

\renewcommand{\arraystretch}{1} %¿ØÖÆ±í¸ñÐÐ¸ßµÄËõ·Å±ÈÀý
\begin{table}[!h]
	\caption{Cost functions from the demand-side}
	%	\fontsize{8}{10}\selectfont
	\begin{tabu} to 1\textwidth {X[1,l] X[10,l] }
		\toprule[1 pt]
		Notation & Description  \\ \midrule
		$\rho_{it}(g_{it})$ & The penalty cost at demand site $i$ in period $t$.\\
        $\rho_{i\infty}(\bm{g}_i)$ & The  penalty cost of eventually unmet demand at demand site $i$,  and $\bm{g}_i = [g_{i1},...,g_{iT}]$.\\
		\bottomrule[1 pt]		
	\end{tabu}
	\label{cost_fun_dem}	
\end{table}

%  \vspace{0.3cm}\begin{tabular}{ll}
%  $f_{jt}$& The facility set-up cost in location $j$ in period $t$.\\
%  $a_{jt}$& The unit cost of expanding capacity in facility  $j$ in period $t$. \\
%  $q_t$& The  limit of supply capacity of every facility in period $t$. \\
%  $c_{ji}$& The unit  cost of transporting relief supplies from location $j$ to demand site $i$.\\
%  % $\beta$& Transportation cost per unit distance.\\
% $\mu(t)$& Unit penalty cost in time period $t$.\\
 % $\alpha$& The exogenous fill rate of the demand satisfaction level, see e.g. type-II service level  \citep{Axsater2015}\\
 %$\eta$& The reliability level of the PD-HL system for all facilities. \\
%  \end{tabular}
 
%  \vspace{0.3cm}\noindent {\bf Random variables and their samples}
 
%  \vspace{0.3cm}\begin{tabular}{ll}
% $d_i$ & The random demand in the demand site $i \in [I]$. \\
% $\hat{d}_i^h$ & An estimate/sample of $d_i$ that was obtained from the $h$-th historical observation.
% \end{tabular}

% \vspace{0.3cm}\noindent {\bf Cost functions from the demand-side}

% \vspace{0.3cm}\begin{tabular}{ll}
% $\rho_{it}(g_{it})$ & The penalty cost at demand site $i$ in period $t$.\\
% $\rho_{i\infty}(\bm{g}_i)$ & The  penalty cost of eventually unmet demand at demand site $i$,  and $\bm{g}_i = [g_{i1},...,g_{iT}]$.
% \end{tabular}

% \vspace{.3cm}

Besides the supply cost of building facilities and transporting relief supplies, we further consider the cost from the demand-side to reflect human sufferings.  Observe that the later the relief supplies are delivered to demand sites, the less benefit of the same amount of supplies from the demand-side point of view. In practice, we should encourage to deliver relief supplies as early as possible. For this purpose, we adopt a monotonically increasing function $\mu(t): \overline{R}\rightarrow R$ to  quantify the timeliness of relief supplies.  Specifically, we consider a logistic growth function \citep{Wang2017poms} that increases sharply at the initial stage and then gradually converges to an upper bound, and propose the following novel penalty functions
\begin{equation*}
% \label{def_Gamma_i}
 \rho_{it} (g_{it})=\mu(t) g_{it}\qquad \text{and}\qquad \rho_{i\infty} (\bm{g}_i)=\mu(\infty)\cdot \left(d_i-\sum\nolimits_{t\in [T]}g_{it}\right)^{+}.
\end{equation*}
%
% . See Figure \ref{fig:mu} for an illustration
%\begin{figure}[htbp]
% \centering
% \includegraphics[width=0.35\linewidth]{fig_mu.eps}
% \caption{A logistic growth function of $\mu(t)$.}
% \label{fig:mu}
%\end{figure} 
%
%
%

\begin{figure}[htbp] \centering
	\subfigure[t = 1] { \label{fig:t1}
		\includegraphics[width=0.47\linewidth]{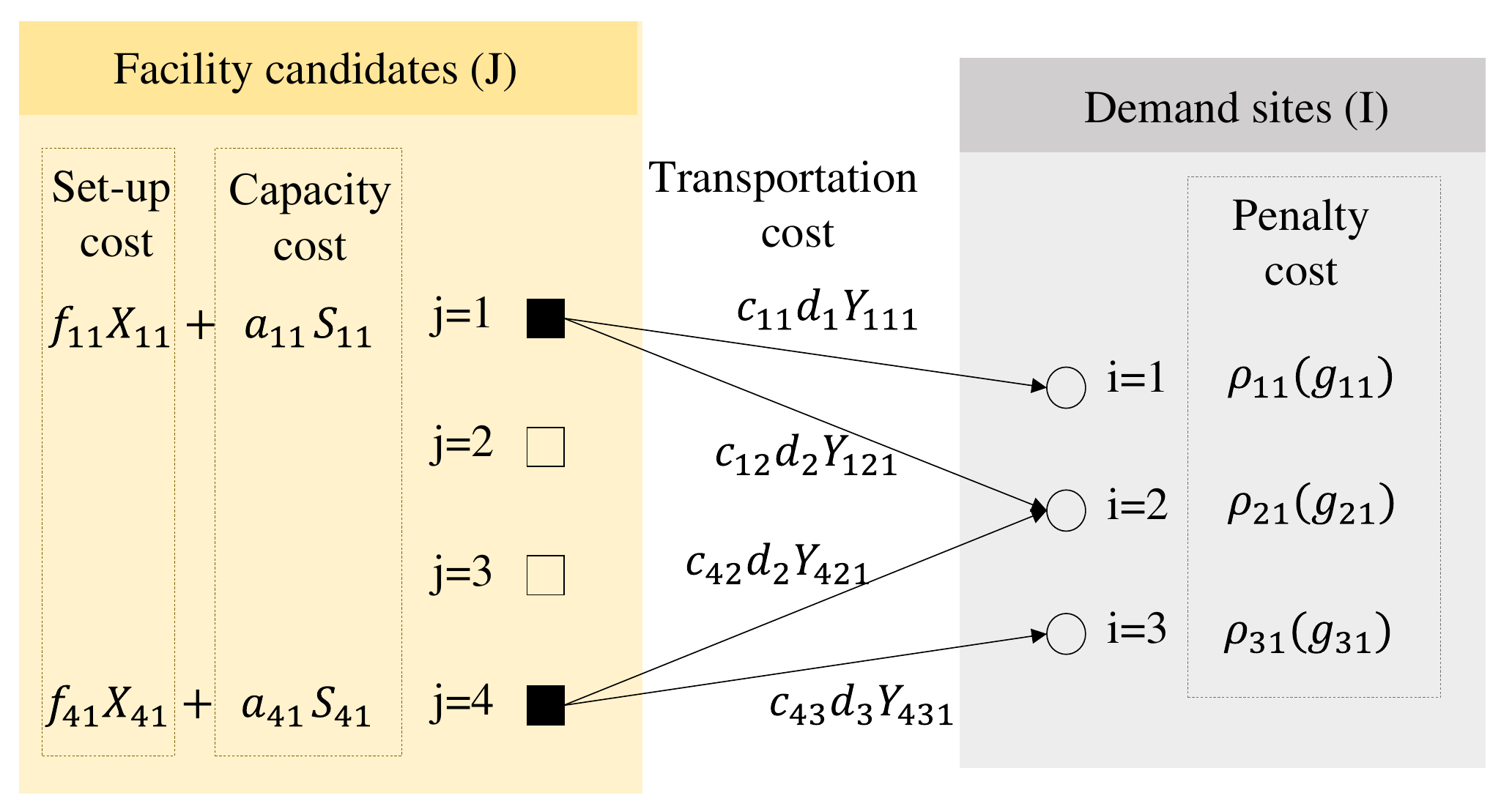}%{cvN=20_k=30.pdf}
	}
	\subfigure[t = 2] { \label{fig:t2}
		\includegraphics[width=0.49\linewidth]{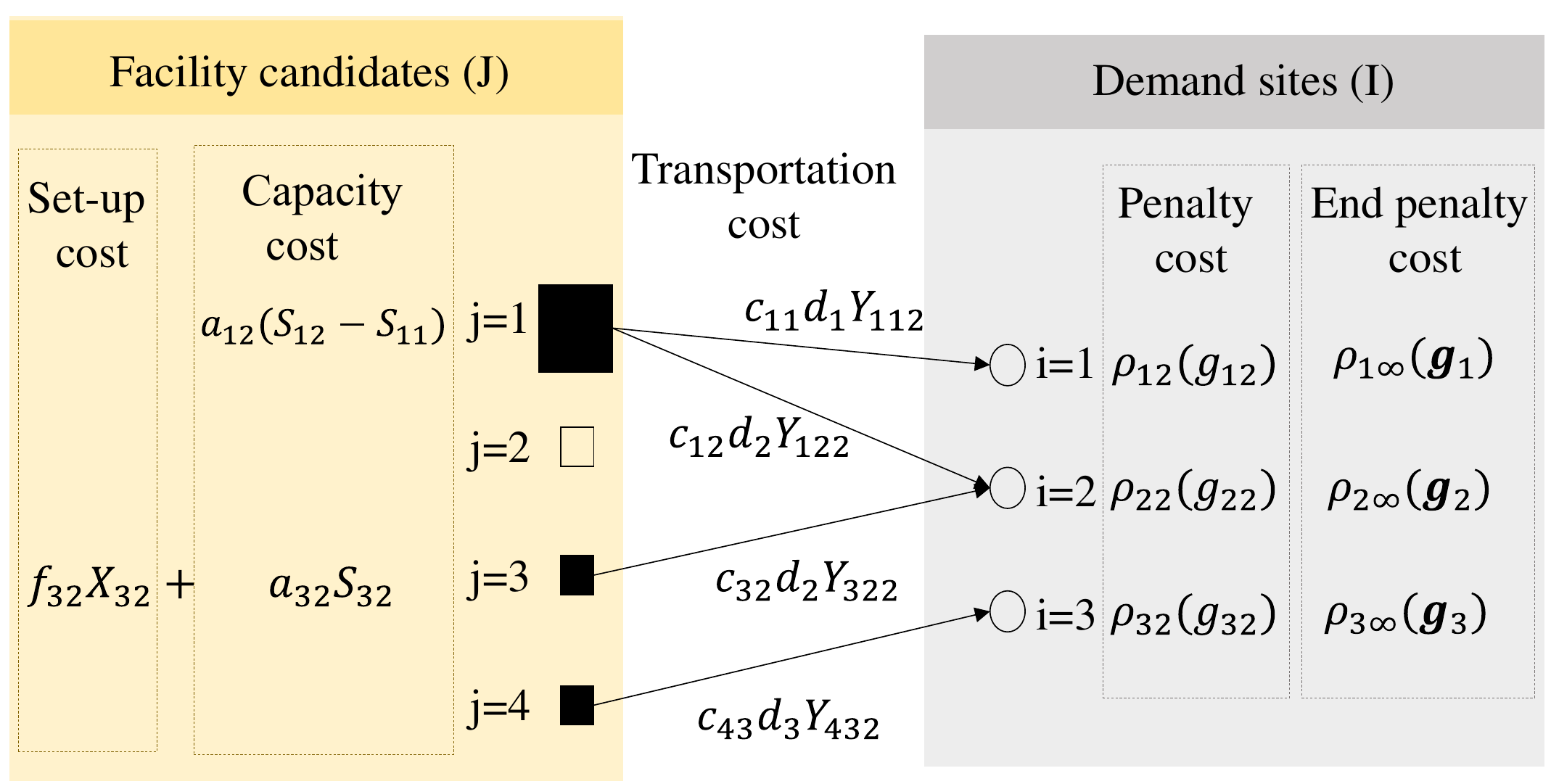}%{cvN=50_k=30.pdf}
	}
	\caption{An example with 4 facility candidates, 3 demand sites and 2 time periods.}
	\label{fig:example432}
\end{figure}

To visualize our objective function,  we provide a toy model with $4$ candidate locations and $3$ demand sites with a $2$-period time horizon  in Figure \ref{fig:example432}.  In period $1$, locations $1$ and $4$ are selected to construct facilities, i.e., $X_{11}=X_{41}=1$, with the associated supply capacities $S_{11}$ and $S_{41}$ respectively, which results in the cost of building facility to be $(f_{11}+a_{11}S_{11})+(f_{41}+a_{41}S_{41})$. While we neglect the cost of getting relief supplies in the facilities\footnote{ To get relief supplies, the facility usually has multiple channels, including self production, donation, purchase from the market, the cost of which are hard to quantify and beyond the scope of this work \citep{CANEL2001411,antunes2001solving,VATSA20211107,testing2021EJOR}. }, another  supply-side cost is from transporting relief supplies, which is given as $(c_{11}d_1Y_{111}+c_{12}d_2Y_{121})+(c_{42}d_2Y_{421}+c_{43}d_3Y_{431})$. We also consider the penalty cost from the demand-side, i.e., $\rho_{11}(g_{11})+\rho_{21}(g_{21})+\rho_{31}(g_{31})$. In period $2$,  location  $3$ is selected to construct facility with supply capacity $S_{32}$ and the capacity of facility $1$ is expanded to $S_{21}$.  Then, the facility cost includes the building cost $(f_{32}+a_{32}S_{32})$ and the expansion cost $a_{12}(S_{21}-S_{11})$.  Similarly, the transportation cost is $(c_{11}d_1Y_{112}+c_{12}d_2Y_{122})+c_{32}d_2Y_{322}+c_{43}d_3Y_{432}$ and the penalty cost from demand-side is $\rho_{12}(g_{12})+\rho_{22}(g_{22})+\rho_{32}(g_{32})$, plus a penalty on the unmet demands, i.e., $\rho_{1\infty}({\bf g}_{1})+\rho_{2\infty}({\bf g}_{2})+\rho_{3\infty}({\bf g}_{3})$. Our  objective aims to minimize the sum of all the above costs under the constraint that $d_i$ is not directly observed. 
 
\subsection{The stochastic MFLCP model with joint chance constraints}
%\addtext{Consider a HL system with $J$ facility candidates $[J] = \{1,2\cdots,J\}$ and $I$ affected areas $[I] = \{1,2\cdots,i\}$. Let $d_i$ denote the demand on each site $i \in [I]$ and $\bm{d} = [d_1,\cdots,d_I]^T$ be the stacked vector of the demand over all sites. }
% where, $\mathbb{F}$ denotes the distribution function of random demand $\bm d$. The joint chance constraint \eqref{joint_chance} ensures the maximum probability of exceeding the capacity over all facility candidates is no larger than a precise level $\eta$. Note that if $\eta = 1$, constraint \eqref{joint_chance} reduces to its deterministic counterpart \eqref{joint_chance} with demand realizations $\bm d$.

%The MFLCP model with joint chance constraints is formulated as follows. First, the objective function \eqref{obj_TC_T} is the summation of the supply-side costs, i.e., total construction cost, capacity expansion cost, transportation cost, and the demand-side costs, i.e., the penalty cost for unmet demand during each period and the final unmet penalty cost at the end of time horizon. 
For a general relief network, the objective function of the MFLCP model  is  given as follows
\begin{equation}
\text{min.}
\sum\limits_{t\in[T]}\left\{\sum\limits_{j\in[J]}f_{jt}X_{jt} +  a_{jt}\left(S_{jt}-S_{j,t-1}\right)+ \mathbb{E}_{\mathbb{F}} \left\{\sum_{i\in[I]}\left(\sum\limits_{j\in[J]} c_{ji}Y_{jit}d_i +\rho_{it}(g_{it})\right)\right\}
\right\} +\mathbb{E}_{\mathbb{F}}\left\{\sum\limits_{i\in[I]}\rho_{i\infty}(\bm{g}_i)\right\}\label{obj_TC_T}
\end{equation}
where the expectation is taken over the random demand vector ${\bf d}=[d_1,\ldots,d_I]^\mathsf{T}$ with unknown distribution $\mathbb{F}$. Next, we describe our constraints of the MFLCP model  in three cases.

 $\bullet$ Trivial constraints.
\begin{equation}
X_{jt}\in\{0,1\}, ~S_{jt} \in \mathbb Z ^+, ~Y_{jit}\ge 0,  \quad\forall  i\in[I], j\in[J],t\in[T] \label{con_X}
\end{equation}
where $Z ^+$ is a set of non-negative integers as we observe that the supply capacity usually takes discrete values. 

 $\bullet$ Reliability constraints
\begin{subequations}
\label{con_reliability}
\begin{align}
&\mathbb{P}_{\mathbb F}\left\{\max_{j \in [J], t \in [T]}\{S_{jt}-\bm{Y}_{jt}^\mathsf{T}\bm{d}\}\ge 0\right\} \ge 1 - \eta, \label{joint_chance}\\
& \alpha\le \sum_{t\in [T]}\sum_{j\in [J]}Y_{jit}\le 1, \forall i\in [I]. \label{con_demand}
\end{align}
\end{subequations}
Intuitively, \eqref{joint_chance} requires that the supply capacity of each facility should be greater than the total of its transported supplies with a high probability $1-\eta$.   \eqref{con_demand} implies that the demand in each site should be served with a percentage  larger than $\alpha$, which is also called Type II service level in \citet{Axsater2015} and less than $1$, which is obviously a trivial constraint. 

 $\bullet$ Capacity constraints
\begin{equation}
\label{SMLCPP}
\sum_{\tau=1}^tX_{j\tau}\le1,~ S_{jt}\le q_t\sum_{\tau=1}^t X_{j\tau},~ S_{jt}\ge S_{j,t-1}, \quad \forall j\in[J],t\in  [T].
\end{equation}
 The first inequality requires that each location can only build at most one facility, the second inequality implies that the supply capacity of  each facility cannot be great than its limit, and the last inequality means that each facility will not reduce its supply capacity for serving relief supplies.

The stochastic MFLCP aims to minimize the objective function \eqref{obj_TC_T} under \eqref{con_X}-\eqref{SMLCPP} and that $\mathbb{F}$ is given.

\subsection{The MFLCP model with $\infty$-Wasserstein joint chance constraints} \label{sec2-3}
Unfortunately, $\mathbb{F}$ in the stochastic MFLCP  is usually unavailable and can only be estimated through the finite historical data $\{\hat{\bm{d}}^h\}_{h\in [H]}$, where $\hat{\bm{d}}^h=[d_1^h,\ldots,d_I^h]^\mathsf{T}$ is an estimate/sample of $\bm{d}$ and $[H] = \{1,\dots,H\}$. Note that these estimates/samples are generally estimated from the population at the demand sites \citep{daskin2011network,zhang2015capacitated}  or collected empirically \citep{blood2014TRE,zhang2015novel}.
Then,  a common idea is to adopt the SAA method to approximate $\mathbb{F}$ by an empirical distribution $\mathbb{F}_H$ over the sample set, i.e.,
$$\mathbb{F}_H({\bm{d}})=\frac{1}{H}\sum_{h\in [H]}\mathbb{I}_{\{\hat{\bm{d}}^{h}= \bm{d}\}},$$
where $\mathbb {I}_A$ is the indicator function of  event $A$. 

%Then, we approximate the MFLCP model \eqref{sto-MFLCP} by
%\begin{equation}
% \label{saa_approx}
% \begin{aligned}
% \text{min.} ~&
%  \sum\limits_{t\in[T]}\left\{\sum\limits_{j\in[J]}f_{jt}X_{jt} +  a_{jt}\left(S_{jt}-S_{j,t-1}\right)+ \frac{1}{H}\sum_{h\in [H]} \left\{\sum\limits_{i\in[I]}\left(\sum\limits_{j\in[J]} \beta c_{ji}Y_{jit}d^h_i +\rho_{it}(g_{it})\right)\right\}\right\} \\ &+\frac{1}{H}\sum_{h\in [H]}\sum_{i\in [I]}\rho_{i\infty}(\bm{g}_i)\\
%  \text{s.t.}&\frac{1}{H}\sum_{h\in [H]}\left\{\max_{j \in [J],t \in [T]}\left\{ \bm{Y}^T_{jt}\bm{\hat{d}}^h - S_{jt} \right\}\le 0\right\}\ge 1- \eta, \\
%  &\eqref{con_demand_T} -\eqref{con_Y}.
% \end{aligned}
%\end{equation}

From the Glivernko-Cantelli Theorem \citep{cantelli1933sulla}, the empirical distribution $\mathbb{F}_H$ converges weakly to the true distribution $\mathbb{F}$ as $H$ increases to infinity, which implies that the SAA method is reliable only if $\mathbb{F}_H$ is a good approximation of $\mathbb{F}$. However, insufficient and/or low-quality samples may result in an unreliable $\mathbb{F}_H$ and a  large deviation from $\mathbb{F}$. Moreover, the demand distribution can not be fixed. In any case, such an SAA approach may exhibit poor out-of-sample performance.

In this work, we adopt a data-driven DR approach to deal with the demand uncertainty. Specifically, we assume that $\mathbb{F}$ belongs to an ambiguity set centered at the empirical distribution $\mathbb{F}_H$ with the radius $\epsilon_H$, and use the $\infty$-Wasserstein metric to evaluate the distance between two distribution functions, leading to the following MFLCP model with $\infty$-Wasserstein joint chance constraints
\begin{subequations}
 \label{dro-p}
 \begin{align}
  \text{MFLCP-W:}\quad\text{min.}\quad&\sum_{t\in[T]}\left\{f_{jt}X_{jt} + \sum_{j\in[J]} a_{jt}\left(S_{jt}-S_{j,t-1}\right)\right\}, \nonumber\\
  &+ \sup_{\mathbb{F}\in\mathcal{F}^{\infty}_H}\mathbb{E}_{\mathbb{F}}\left\{\sum_{t \in [T]}\sum_{i \in [I]}\left(\sum_{j\in [J]} c_{ji}Y_{jit}d_i+\rho_{it}(g_{it})\right)+\sum_{i \in [I]}\rho_{i\infty}(\bm{g}_i)\right\}\label{obj_dro} \\
  \text{s.t.}\quad&\inf_{\mathbb{F}\in\mathcal{F}^{\infty}_H}\mathbb{P}_{\mathbb{F}}\left\{\max_{j \in [J],t \in [T]}\left\{ \bm{Y}_{jt}^\mathsf{T}\bm d - S_{jt} \right\}\le 0\right\}\ge 1 - \eta, \label{con_joint_dro}\\
  \quad&\eqref{con_X}, \eqref{con_demand}, \eqref{SMLCPP} \label{Trival_Con}.
 \end{align}
\end{subequations}

\subsection{The ambiguity set via $\infty$-Wasserstein metric}
We introduce the following metric for the construction of the Wasserstein ball $\mathcal{F}_H^\infty$.

\begin{definition}\citep{kantorovich1958space}
Let $dis(\bm \xi^1,\bm \xi^2)=||\bm \xi^1-\bm \xi^2||_p$, $(\Xi,dis)$ be a Polish metric space and $\mathcal{M}(\Xi)$ denote the support set of some probability distributions. The type $r$-Wasserstein distance between $\mathbb F_1\in \mathcal M (\Xi)$ and $\mathbb F_2\in \mathcal M (\Xi)$ is defined as
 \begin{equation*}
  W^r(\mathbb F _1,\mathbb F_2)=\inf_{\mathbb{K}\in \mathcal{K}(\mathbb{F}_1,\mathbb{F}_2)} \left\{\left(\int_{\Xi\times \Xi}dis(\bm \xi^1,\bm \xi ^2)^r \mathbb{K}(d\bm \xi^1,d \bm \xi^2)\right)^{\frac{1}{r}}: \int_\Xi \mathbb{K}(\bm \xi^1,d\bm \xi^2)=\mathbb F_1(\bm\xi^1),\int_\Xi \mathbb{K}(d\bm\xi^1,\bm \xi^2)=\mathbb F_2(\bm\xi^2)\right\},
 \end{equation*}
and the $\infty$-Wasserstein metric is the limit of $W^r$ as $r$ tends to infinity and amounts to
 \begin{equation}
 \label{W_infty}
  W^\infty(\mathbb F_1,\mathbb F_2)=\inf \left\{\esssup_{\mathbb{K}}\left\{dis(\bm \xi^1, \bm \xi^2):(\bm \xi^1,\bm \xi^2)\in \Xi \times \Xi\right\}: \mathbb{K}\in \mathcal{K}(\mathbb{F}_1,\mathbb{F}_2) \right\},
 \end{equation}
 \noindent where $\mathcal{K}(\Xi\times\Xi)$ is the set of the joint distributions on $\Xi\times\Xi$ with marginal distribution $\mathbb F_1$ and $\mathbb F_2$ respectively. 
\end{definition}

Motivated by recent works on DR optimization with $\infty$-Wasserstein set \citep{chen2018data,wang2021wasserstein}, we adopt the $\infty$-Wasserstein metric and $2$-norm, i.e., $r = \infty$ and $p = 2$ in $dis(\cdot,\cdot)$ to construct the Wasserstein ball $\mathcal{F}_H^\infty$, i.e.,
\begin{equation}
\begin{aligned}
\label{Wset}
 &\mathcal{F}^{\infty}_H=
 &\left\{\begin{array}{l|l}\mathbb{F} \in \mathcal{M}(\Xi) & \begin{array}{ll}\bm{d} \sim \mathbb{F}, W^\infty(\mathbb{F},\mathbb{F}_{H})\le \epsilon_H \\
 \end{array}\end{array}\right\}.
\end{aligned}
\end{equation}
Obviously, the ambiguity set $\mathcal{F}^{\infty}_H$ consists of all distributions within $\epsilon_H$-distance from the empirical distribution $\mathbb{F}_H$. The radius $\epsilon_H$ reflects the confidence level on $\mathbb{F}_H$, and the sample dataset can be utilized more flexibly to control the conservatism of the distribution uncertainty via tuning the radius.

The MFLCP-W model is an infinite-dimensional problem and appears to be intractable. However, in sharp contrast to the infinite-dimensional DR joint chance constrained models in \citet{xie2021distributionally,chen2018data}, which  are reformulated as intractable bilinear problems and can only be solved by providing feasible solutions. Our model can be equivalently reformulated as a finite-dimensional MISOCP. 

However, even the basic location models which are formulated as MIP problems are NP-Hard \citep{lewis1979computers,current2002discrete}, our MISOCP model is far more complex and difficult to be directly solved by commercial solvers for large-scaled instances. To this end, we design an OA algorithm with two-tailored cuts to efficiently solve our MISOCP especially for the large-sized instances. Our algorithm is able to obtain the optimal solution within finite number of iterations, as the continuous relaxation of MISOCP is a convex program \citep{BONAMI2008186}.

To the best of our knowledge, we are the first to adopt the $\infty$-Wasserstein joint chance constraints in the MFLCP model. The most relevant paper that involves the Wasserstein set in a facility location problem is \cite{saif2021data}. They propose a single-period model and require solution to be feasible for all demand realizations, which is more conservative than ours with joint chance constraints.

\section{Reformulation of the MFLCP-W model}\label{section:reformulation}
In this section, we reformulate the infinite-dimensional MFLCP-W model \eqref{dro-p} over the Wasserstein ball $\mathcal{F}_H^\infty$ in \eqref{Wset} to a finite MISOCP. Before the reformulation, we express the objective function in \eqref{dro-p} as a vector form at first.

\begin{lemma}
The objective function \eqref{obj_dro} can be expressed in a vector form as follows:
\begin{equation}
\label{obj_dro_vector}
 \min. \quad \sum_{t\in[T]}\left\{f_{jt}X_{jt} + \sum_{j\in[J]} a_{jt}\left(S_{jt}-S_{j,t-1}\right)\right\} +\sup_{\mathbb{F}\in\mathcal{F}^{\infty}_H}\mathbb{E}_{\mathbb{F}}\left\{\bm u^\mathsf{T}\bm{d} +\bm v^\mathsf{T}\bm{d}+\bm{z}^\mathsf{T}\bm{d} \right\}
\end{equation}
where
\begin{subequations}
 \begin{align}
 & u_{i} = \sum_{t \in [T]}\sum_{j\in [J]} c_{ji}Y_{jit}, ~~ \forall i \in [I] \label{con_u}\\
 & v_{i} = \sum_{t \in [T]}\sum_{j\in [J]}\mu(t)Y_{jit}, ~~ \forall i \in [I] \label{con_v}\\
 & z_{i} = \mu(\infty)\left(1-\sum_{t \in [T]}\sum_{j \in [J]}Y_{jit}\right), ~~ \forall i \in [I] \label{con_z}\\
 & u_i,v_i,z_i\in\mathbb R,\forall i\in[I].\label{con_uvz}
 \end{align}
\end{subequations}
Moreover, the MFLCP-W model \eqref{dro-p} is equivalent to
\begin{align}
  \min.~& \sum_{t\in[T]}\left\{f_{jt}X_{jt} + \sum_{j\in[J]} a_{jt}\left(S_{jt}-S_{j,t-1}\right)\right\} +\sup_{\mathbb{F}\in\mathcal{F}^{\infty}_H}\mathbb{E}_{\mathbb{F}}\left\{\bm u^\mathsf{T}\bm{d} +\bm v^\mathsf{T}\bm{d}+\bm{z}^\mathsf{T}\bm{d} \right\} \label{dro-pm} \\
 \sta ~& \eqref{con_joint_dro}, \eqref{Trival_Con},  \eqref{con_u}-\eqref{con_uvz}.\nonumber
\end{align}
\end{lemma}

\begin{proof}
Observe that
\begin{equation*}
 \begin{aligned}
 \sum_{t\in [T]}\rho_{it}\left(g_{it}\right)=\sum_{t\in [T]}\mu(t)g_{it}=\sum_{t\in [T]}\mu(t)\left(\sum_{j\in [J]}Y_{jit}d_i\right)=\left(\sum_{t\in [T]}\sum_{j\in [J]}\mu(t)Y_{jit}\right)d_i
 \end{aligned}
\end{equation*}
and
\begin{equation*}
 \begin{aligned}
 \rho_{i\infty}\left(\bm{g}_i\right)=\mu(\infty)\max\left\{d_i-\sum_{t\in [T]}g_{it},0\right\}=\mu(\infty)\max\left\{d_i-\sum_{t\in [T]}\sum_{j\in [J]}Y_{jit}d_i,0\right\}=\mu(\infty)\left(1-\sum_{t\in [T]}\sum_{j\in [J]}Y_{jit}\right)d_i,
 \end{aligned}
\end{equation*}
where the last equality holds due to constraint \eqref{con_demand}.

Thus, the expectation  in the objective function of \eqref{obj_dro} can be expressed as
\begin{equation*}
 \begin{aligned}
 &\mathbb{E}_{\mathbb{F}}\left\{\sum_{t \in [T]}\sum_{i \in [I]}\left(\sum_{j\in [J]}c_{ji}Y_{jit}d_i+\rho_{it}(g_{it})\right)+\sum_{i \in [I]}\rho_{i\infty}(\bm{g}_i)\right\} \\
 &=\mathbb{E}_{\mathbb{F}}\left\{\sum_{t \in [T]}\sum_{i \in [I]}\sum_{j\in [J]} c_{ji}Y_{jit}d_i+ \sum_{t \in [T]}\sum_{i \in [I]}\rho_{it}(g_{it})+\sum_{i \in [I]}\rho_{i\infty}(\bm{g}_i)\right\} \\
 &=\mathbb{E}_{\mathbb{F}}\left\{\sum_{i \in [I]}\sum_{t \in [T]}\sum_{j\in [J]} c_{ji}Y_{jit}d_i+ \sum_{i \in [I]}\sum_{t\in [T]}\sum_{j\in [J]}\mu(t)Y_{jit}d_i+\sum_{i \in [I]}\mu(\infty)\left(1-\sum_{t\in [T]}\sum_{j\in [J]}Y_{jit}\right)d_i\right\}\\
 &=\mathbb{E}_{\mathbb{F}}\left\{\sum_{i \in [I]}(u_id_i+v_id_i+z_id_i)\right\}=\mathbb{E}_{\mathbb{F}}\{\bm u^\mathsf{T}\bm d +\bm v^\mathsf{T}\bm d +\bm z^\mathsf{T}\bm d\},
 \end{aligned}
\end{equation*}
where $u_i = \sum_{t\in [T]}\sum_{j \in [J]}c_{ji}Y_{jit}$, $v_i = \sum_{t\in [T]}\sum_{j \in [J]}\mu(t)Y_{jit}$ and $z_i = \mu(\infty)\left\{1-\sum_{t\in [T]}\sum_{j \in [J]}Y_{jit}\right\}$.
Thus, replacing \eqref{obj_dro} with \eqref{obj_dro_vector} and adding the corresponding constraints \eqref{con_u}-\eqref{con_uvz} to \eqref{dro-p} result in an equivalent problem \eqref{dro-pm}.
\end{proof}

Then, we show that model \eqref{dro-pm} over the $\infty$-Wasserstein ambiguity set $\mathcal{F}_H^\infty$ defined in \eqref{Wset} can be reformulated as an MISOCP.

\begin{theorem}
\label{theo_equi}
The reformulated model \eqref{dro-pm} over the $\infty$-Wasserstein ball $\mathcal{F}_H^\infty$ is equivalent to an MISOCP  as follows,
\begin{subequations}
\label{dro-equi-socp}
 \begin{align}
  \min.\quad &\sum_{t\in[T]}\left\{f_{jt}X_{jt} + \sum_{j\in[J]} a_{jt}\left(S_{jt}-S_{j,t-1}\right)\right\}+\frac{1}{H}\sum_{h \in [H]}l_h, \\
  \text{\rm{s.t.}}\quad &\epsilon_H\|\bm u+\bm v+\bm{z}\|_2+(\bm u+\bm v+\bm{z})^\mathsf{T}\hat{\bm{d}}^h\le l_h, \forall h \in [H],\label{dro_cone1}\\
  &\epsilon_H\|\bm{Y}_{jt}\|_2+\bm{Y}_{jt}^\mathsf{T}\hat{\bm{d}}^h-S_{jt}\le M_h(1-w_h), ~\forall j \in [J], t \in [T], h \in [H], \label{dro_cone2} \\
   &\frac{1}{H}\sum_{h \in [H]}w_h \ge 1-\eta,\label{dro_theta} \\
  &w_h \in \{0,1\}, ~\forall h \in [H], \label{dro_theta_domain}\\
  &\eqref{Trival_Con}, \eqref{con_u}-\eqref{con_uvz},\nonumber
 \end{align}
\label{dro_dflp_redef}
\end{subequations}
where $M_h$ is the big-M constant.
\end{theorem}

Before the proof for Theorem \ref{theo_equi}, a lemma is introduced below.

\begin{lemma}
\label{lemma_inf_equi}
Fix a collection of samples $\{\hat{\bm{d}}^h\}_{h \in [H]}$, a Wasserstein radius $\epsilon_H > 0$, a transport plan $\bm{Y}$ and capacity plan $\bm{S}$, $\frac{1}{H}\sum_{h \in [H]}\inf_{\|\bm{d}-\hat{\bm{d}}^h\|_2 \le \epsilon_H}\mathbb{I}\left\{\max_{j\in[J], t\in[T]}\left\{\epsilon_H\|\bm{Y}_{jt}\|_2+\bm{Y}_{jt}^\mathsf{T}\bm{d}-S_{jt}\right\}\le 0\right\}$ is equivalent to the maximum value of the following mixed 0-1 program,
\begin{equation}
 \begin{aligned}
 \label{chance_equi_MLP}
 {\max}. \quad &\frac{1}{H}\sum_{h \in [H]}w_h, \\
 \sta \quad &\epsilon_H\|\bm{Y}_{jt}\|_2+\bm{Y}_{jt}^\mathsf{T}\hat{\bm{d}}^h-S_{jt}\le M_h(1-w_h), \\
 &w_h \in \{0,1\}, ~\forall j \in [J], t \in [T], h \in [H].
 \end{aligned}
\end{equation}
% where $M_h$ is the big-M constant.
% $||\cdot||_q$ denotes the dual norm of $||\cdot||_p$, and $\frac{1}{p}+\frac{1}{q}=1$ and
\end{lemma}

\begin{proof}
With fixed sample $\hat{\bm{d}}^h$, for any $\bm{Y}_{jt}$ and $S_{jt}$, we have
\begin{equation}\label{def_indicator}
\begin{aligned}
	&\inf_{\|\bm{d}-\hat{\bm{d}}^h\|_2\le \epsilon_H}\mathbb{I}\left\{\max_{j \in [J], t \in [T]}\{\bm{Y}_{jt}^\mathsf{T}\bm{d}-S_{jt}\}\le 0\right\},\\
	=&\left\{
		\begin{array}{ll}
		1 & \text{if} ~ \bm{Y}_{jt}^\mathsf{T}\bm{d}\le S_{jt}, \forall \ \|\bm{d}-\hat{\bm{d}}^h\|_2\le \epsilon_H , j \in [J], t \in [T]\\
		0 & \text{otherwise}
		\end{array}
	\right.,\\
	=&\left\{
		\begin{array}{ll}
		1 & \text{if} ~ \sup_{\|\bm{d}-\hat{\bm{d}}^h\|_2\le \epsilon_H}\bm{Y}_{jt}^\mathsf{T}\bm{d}\le S_{jt}, \forall j \in [J], t \in [T]\\
		0 & \text{otherwise}
		\end{array}
	\right.,\\
	=&\left\{
		\begin{array}{ll}
		1 & \text{if} ~ \epsilon_H\|\bm{Y}_{jt}\|_2+\bm{Y}_{jt}^\mathsf{T}\hat{\bm{d}}^h \le S_{jt}, \forall j \in [J], t \in [T]\\
		0 & \text{otherwise}
	\end{array}
	\right.,\\
\end{aligned}
\end{equation}
where the last equality follows from the definition of dual norm,
\begin{equation}
\begin{aligned}
\sup_{\|\bm{d}-\hat{\bm{d}}^h\|_2\le \epsilon_H}\bm Y_{jt}^\mathsf{T}\bm{d}=&\sup_{\|\Delta\|_2\le \epsilon_H}\bm Y_{jt}^\mathsf{T}(\Delta + \hat{\bm{d}}^h)
=\sup_{\|\Delta\|_2\le \epsilon_H}\bm Y_{jt}^\mathsf{T}\Delta + \bm Y_{jt}^\mathsf{T}\hat{\bm{d}}^h \\
=&\epsilon_H\|\bm Y_{jt}\|_2+\bm Y_{jt}^\mathsf{T}\hat{\bm{d}}^h.
\end{aligned}
\end{equation}
Then, the equivalent form of \eqref{chance_equi_MLP} can be obtained through the epigraphical formulation of each infimum term. Specifically, the equivalence follows that at optimality, the $w_h = 1$ if the infimum of the indicator function in \eqref{def_indicator} equals to $1$ and $w_h = 0$ otherwise.
\end{proof}

Next we provide the proof of Theorem \ref{theo_equi}.

\begin{proof} [Proof of Theorem \ref{theo_equi}]
We first reformulate the objective function \eqref{obj_dro_vector} as a second-order conic program (SOCP) for any fixed $\bm{u},\bm{v}$ and $\bm{z}$.
According to the Theorem 5 of \citet{bertsimas2018data},
\begin{equation*}
 \begin{aligned}
 &\sup_{\mathbb{F}\in\mathcal{F}^\infty_H}\mathbb{E}_{\mathbb{F}}\left[\bm u^\mathsf{T}\bm{d} +\bm v^\mathsf{T}\bm{d}+\bm{z}^\mathsf{T}\bm{d} \right]
 =\frac{1}{H}\sum_{h \in [H]}\sup_{\|\bm{d}-\hat{\bm{d}}^h\|_2\le \epsilon_H}\left\{\bm u^\mathsf{T}\bm{d} +\bm v^\mathsf{T}\bm{d}+\bm{z}^\mathsf{T}\bm{d} \right\}.
 \end{aligned}
\end{equation*}
By the properties of dual norm, the supremum term of each sample $\hat {\bm{d}}^h$ is equivalent to
\begin{equation*}
 \sup_{\|\bm{d}-\hat{\bm{d}}^h\|_2\le \epsilon_H}\left\{\bm u^\mathsf{T}\bm{d} +\bm v^\mathsf{T}\bm{d}+\bm{z}^\mathsf{T}\bm{d}\right\} = \epsilon_H\|\bm u+\bm v+\bm{z}\|_2+(\bm u+\bm v+\bm{z})^\mathsf{T}\hat{\bm{d}}^h.
\end{equation*}

From the epigraphical reformulation, the objective function in \eqref{dro-pm} is reformulated as
\begin{equation}
 \begin{aligned}
 \label{euqi-obj}
 \text{min.}\quad &\sum_{t\in[T]}\left\{f_{jt}X_{jt} + \sum_{j\in[J]} a_{jt}\left(S_{jt}-S_{j,t-1}\right)\right\}+\frac{1}{H}\sum_{h \in [H]}l_h,\\
 \text{s.t.}\quad &\epsilon_H\|\bm u+\bm v+\bm{z}\|_2+(\bm u+\bm v+\bm{z})^\mathsf{T}\hat{\bm{d}}^h\le l_h,\forall h\in[H]. \\
 \end{aligned}
\end{equation}
 
Next, we provide the reformulation of joint chance constraint \eqref{con_joint_dro}.
Let 1 minus both sides of the inequality \eqref{con_joint_dro}, we obtain that
$$
\inf_{\mathbb{F}\in \mathcal{F}_H^\infty}\mathbb{P}_{\mathbb{F}}\left\{\max_{j \in [J],t \in [T]}\left\{ \bm{Y}_{jt}^\mathsf{T}\bm{d} - S_{jt} \right\}\le 0\right\} \ge 1-\eta \iff \sup_{\mathbb{F}\in \mathcal{F}_H^\infty}\mathbb{P}_{\mathbb{F}}\left\{\max_{j \in [J],t \in [T]}\left\{ \bm{Y}_{jt}^\mathsf{T}\bm{d} - S_{jt} \right\} > 0\right\} \le \eta.
$$
Note that we have
\begin{equation*}
    \mathbb{P}_{\mathbb{F}}\left\{\max_{j \in [J],t \in [T]}\left\{ \bm{Y}_{jt}^\mathsf{T}\bm{d} - S_{jt} \right\}> 0\right\} = \mathbb{E}_{\mathbb{F}}\left\{\mathbb{I}\left\{\max_{j \in [J], t \in [T]}\{\bm{Y}_{jt}^\mathsf{T}\bm{d}-S_{jt}\}> 0\right\}\right\}.
\end{equation*}
Then by applying Theorem 5 of \citet{bertsimas2018data} again, the constraint \eqref{con_joint_dro} is equivalent to
\begin{equation}
 \label{equi_chance_m}
 \frac{1}{H}\sum_{h \in [H]}\sup_{\|\bm{d}-\hat{\bm{d}}^h\|_2\le \epsilon_H}\mathbb{I}\left\{\max_{j \in [J], t \in [T]}\{\bm{Y}_{jt}^\mathsf{T}\bm{d}-S_{jt}\} > 0\right\} \le \eta.
\end{equation}
Using the fact that $\mathbb{I}\left\{\max_{j \in [J], t \in [T]}\{\bm{Y}_{jt}^\mathsf{T}\bm{d}-S_{jt}\}\le 0\right\} + \mathbb{I}\left\{\max_{j \in [J], t \in [T]}\{\bm{Y}_{jt}^\mathsf{T}\bm{d}-S_{jt}\}> 0\right\} = 1$, the \eqref{equi_chance_m} is equivalent to 
\begin{equation}
\label{equi_chance_1}
\frac{1}{H}\sum_{h \in [H]}\inf_{\|\bm{d}-\hat{\bm{d}}^h\|_2\le \epsilon_H}\mathbb{I}\left\{\max_{j \in [J], t \in [T]}\{\bm{Y}_{jt}^\mathsf{T}\bm{d}-S_{jt}\}\le 0\right\} \ge 1-\eta.
\end{equation}
%$$\frac{1}{H}\sum_{h=1}^{H}\sup_{\|\bm{d}-\hat{\bm{d}}^h\|_2\le \epsilon_H}\mathbb{I}\left\{\max_{j \in [J], t \in [T]}\{\bm{Y}_{jt}^\mathsf{T}\bm{d}-S_{jt}\}> 0\right\} \le \eta,$$
%which is equivalent to

According to Lemma \ref{lemma_inf_equi}, setting the
optimal value of the problem in \eqref{chance_equi_MLP} to be no less than $(1- \eta)$ leads to the equivalence between the constraint \eqref{equi_chance_1} and the following constraints
\begin{equation}
\begin{aligned}
\label{equi_chance_2}
~&\frac{1}{H}\sum_{h \in [H]}w_h \ge 1-\eta, \\
~&\epsilon_H\|\bm{Y}_{jt}\|_2+\bm{Y}_{jt}^\mathsf{T}\hat{\bm{d}}^h-S_{jt}\le M_h(1-w_h), ~\forall j \in [J], t \in [T], h \in [H],  \\
&~w_h \in \{0,1\}, ~ \forall h \in [H].
\end{aligned}
\end{equation}
Consequently, we reformulate the DR joint chance constraint \eqref{con_joint_dro} as second-order conic constraints with big-M constants.

By substituting the equivalent formulation \eqref{equi_chance_2} of the $\infty$-Wasserstein joint chance constraints \eqref{con_joint_dro} and the equivalent formulation \eqref{euqi-obj} of the objective function into \eqref{dro-pm}, we reformulate the MFLCP-W model \eqref{dro-pm} over the $\infty$-Wasserstein ball $\mathcal{F}_H^\infty$ as the MISOCP \eqref{dro-equi-socp}.
\end{proof}

We further derive a lower bound of the big-M coefficients $M_h$ by inspection. For example, since $\bm Y_{jt} \in [0,1]^I$ and $S_{jt} \ge 0$, then for each $h \in [H]$, one possible $M_h$ can be derived as bellow:
\begin{equation*}
\begin{aligned}
M_h = \max_{j \in [J], t \in [T]}\left\{\epsilon_H\|\bm{e}_I\|_2 + \sum_{i\in[I]}d^h_i\right\},
\end{aligned}
\end{equation*}
where $\bm{e}_I \in \mathbb{R}^I$ is the all one vector. To tighten the big-M coefficients, interested readers could refer to \citet{qiu2014covering} and \citet{song2014chance} for more details.

Theorem \ref{theo_equi} implies that the MFLCP-W model \eqref{dro-pm} over the $\infty$-Wasserstein ball is equivalent to a finite MISOCP. Although problem \eqref{dro-equi-socp} can be solved directly by the commercial solvers, it is time-consuming for large-sized instances. Thus,  we design an OA algorithm with two tailored cuts to efficiently solve problem \eqref{dro-equi-socp} in the following. 

%Furthermore, different $p$-norms in Wasserstein distance lead to different equivalent problems of our model. For example, if $p=1$ in $dis(\cdot,\cdot)$, the proposed model is equivalent to a mixed integer linear program (MILP).

\section{Solution approach}\label{section:oa}
%\highlight{add the necessity for using OA}
The OA algorithm \citep{1987An} is efficient for solving mixed integer nonlinear programs (MINLP). As the MFLCP-W model is equivalent to an MISOCP, a special case of the MINLP, we utilize the OA algorithm to solve this problem. Moreover, we enhance the OA method by two problem-specific cuts, which improve the computational efficiency. Particularly, we solve the problem \eqref{dro_dflp_redef} by a master-subproblem framework, that is, decomposing the original MISOCP into a series of MILP (named as master problem, MP) and nonlinear programs (named as subproblem, SP) and obtain the optimal solution by iteratively solving SP and MP.

OA algorithm can obtain an optimal solution within finite number of iterations if the continuous relaxation of MINLP is proved to be convex \citep{BONAMI2008186}. Thus, before introducing the detailed procedure of OA, we firstly prove that the continuous relaxation of problem \eqref{dro_dflp_redef} is convex. Note that, constraints \eqref{dro_cone1} and \eqref{dro_cone2} are the only nonlinear components of the reformulated DR model, thus, Lemmas \ref{lemma:convex1} and \ref{lemma:convex2} focus on the convexity of  \eqref{dro_cone1} and continuous relaxations of \eqref{dro_cone2}, i.e., $\phi_h(\bm u,\bm v,\bm z,l_h)$ and $\psi_{jth}(\bm Y_{jt},w_h)$, respectively.

\begin{lemma}
\label{lemma:convex1}
If $\bm u\in\mathbb R^I$, $\bm v\in\mathbb R^I$, $\bm z\in\mathbb R^I$, $l_h\in \mathbb R$, the function
\begin{equation}
 \label{def_phi}
 \phi_h(\bm u,\bm v,\bm z,l_h)=\epsilon_H||\bm u+\bm v+\bm z||_2+(\bm u+\bm v+\bm z)^\mathsf{T}\hat{\bm{d}}^h-l_h
\end{equation}
is convex for any $h \in [H]$.
\end{lemma}
\begin{proof}
For the ease of exposition, we omit index $h$ in the proof. Define an auxiliary decision variable $\bm A\in \mathbb R^I$ such that $\bm A=\bm u+\bm v+\bm z$, then, $\phi(\bm u,\bm v,\bm z,l)$ can be rewritten as $\phi(\bm A,l)=||\bm A||_2+\bm A^\mathsf{T}\hat{\bm d}-l$. As a result, the  convexity map transforms to prove the only nonlinear term $||\bm A||_2$ is convex. Let $||\cdot||_2$ be the Euclidean norm, according to the triangle inequality of norm functions, for any feasible vectors $\bm A_1,\bm A_2\in \mathbb R^I$ and $\lambda\in[0,1]$, it is obvious that, $||\lambda\bm A_1+(1-\lambda)\bm A_2||_2\le \lambda||\bm A_1||_2+(1-\lambda)||\bm A_2||_2,$ thus, convexity holds by definition.
\end{proof}

\begin{lemma}
\label{lemma:convex2}
If $\bm Y_{jt}\in [0,1]^I$, $w_h\in[0,1]$, the function
\begin{equation}
 \label{def_psi}
 \psi_{jth}(\bm Y_{jt},w_h)=\epsilon_H||\bm Y_{jt}||_2+\bm Y_{jt}^\mathsf{T}\hat{\bm{d}}^h-S_{jt}-M_h(1-w_h)
\end{equation}
is convex for any $j\in[J], t\in[T],h\in[H]$.
\end{lemma}
\begin{proof}
The linearity of convexity suffices to show the first term $||\bm Y||_2$ is convex, and is similar to the proof of Lemma \ref{lemma:convex1}.
\end{proof}

Then, we address the details of our tailored OA algorithm in the aspects of initialization, MP and SP from Sections \ref{subsec:init} to \ref{subsec:MP}. Finally, the complete OA procedure is summarized in Section \ref{subsec:procedure}.

Let $K$ denote the maximum iteration number of the OA algorithm, and $k\in[K]=\{0,1,\cdots, K\}$ be its index. %Variables indexed by $0$ represents its initial solution. In addition, define $\bm e$ as a vector where all its elements equal to 1.
The decision variables labeled with a hat indicate the values obtained along with the iterations, e.g., $\hat{\bm S}^k$ represents the integer values obtained by the $k^{th}$ iteration of MP.

%\begin{equation}
% \begin{aligned}
%  \text{min.}&\sum_{t\in[T]}f_{jt}X_{jt} + \sum_{t\in[T]}\sum_{j\in[J]} %a_{jt}\left(S_{jt}-S_{j,t-1}\right)+\beta\frac{1}{H}\sum_{h=1}^{H}\eta_h \\
%  \text{s.t.} ~ &\epsilon_H\|\bm u+\bm v+\bm{z}\|_q+(\bm u+\bm v+\bm{z})^\mathsf{T}\bm{d}_h\le \eta_h, \forall h \in [H] \\
%  &\epsilon_H\|\bm{Y}_{jt}\|_q+\bm{Y}_{jt}^\mathsf{T}\bm{d}_h-S_{jt}\le M_h(1-\theta_h), ~\forall j \in [J], t \in [T], h \in [H], \\
%   &\frac{1}{H}\sum_{h=1}^{H}\theta_h \ge 1-\eta \\
%& \sum_{t\in [T]}\sum_{j\in \text{Cov}_i}Y_{jit}\ge \alpha,\quad \forall i\in[I],\\
%%&\sum_{t\in [T]}\sum_{j\in \text{Cov}_i}Y_{jit}\le 1,\quad \forall i\in[I],\\
%&S_{jt}\le q_t\sum_{\tau=1}^\mathsf{T} X_{j\tau},\quad \forall j\in[J], t\in [T],\\
%&\sum_{\tau=1}^tX_{j\tau}\le1, \forall j\in[J],t\in\mathcal T,\\
%& u_{i} = \sum_{t \in [T]}\sum_{j\in J_i}c_{ji}Y_{jit}, ~~ \forall i \in [I] \\
% & v_{i} = \sum_{t \in [T]}\sum_{j\in J_i}\mu(t)Y_{jit}, ~~ \forall i \in [I] \\
% & z_{i} = \mu(\infty)\left\{1-\sum_{t \in [T]}\sum_{j \in J_i}Y_{jit}\right\}, ~~ \forall i \in [I]\\
% &S_{jt}\ge S_{j,t-1},\quad \forall j\in[J],t\in[T],\\
%&X_{j}\in\{0,1\}, \quad\forall j\in[J],t\in[T],\\
%& S_{jt} \ge 0, \quad \forall j\in[J],t\in[T],\\
%&Y_{jit}\ge 0, \quad \forall i\in[I], j\in[J],t\in[T].\\
%&\theta_h \in \{0,1\}, h \in [H], \\
%\end{aligned}
%\end{equation}

\subsection{Initialization} \label{subsec:init}
We myopically introduce a conservative capacity planning strategy as the initial solution. Let the capacity in each period ($S_{jt}$) equals to the maximum amount that can be established ($q_t$), and then the joint chance constraints can be satisfied simultaneously with probability 1.
i.e.,
$$\hat S_{jt}^0=q_t,\forall j\in[J], t\in [T],\qquad \hat w^0_h=1,\quad\forall h\in[H].$$
By definition, $\hat X_{jt}^0$ can be directly reflected by the values of $\hat S_{jt}^0$, i.e.,
\begin{eqnarray*}\hat X_{jt}^0=\left\{
\begin{array}{cc}
  1,& \text{if } ~~ t = 1 ~~ \\
  0,& \text{otherwise }
\end{array}\right.,\forall j\in[J],t\in[T].
\end{eqnarray*}
\subsection{Subproblem}\label{subsec:SP}
SP finds the optimal solution of the continuous variables with fixed integer input, and provides an upper bound (UB) to the problem \eqref{dro_dflp_redef}. Let $(\hat X_{jt}^k, \hat{S}^k_{jt},\hat{w}^k_h)$ be the integer values obtained in the $k^{th}$ iteration of MP, SP admits the optimal value of continuous variables $\bm Y, \bm u, \bm v, \bm z$ and $\bm l$. The second-order conic SP in the $k^{th}$ iteration can be summarized as follows.
\begin{subequations}
\begin{align}
\quad \text{min.}\quad & \frac{1}{H}\sum_{h=1}^Hl_h,\label{sp_obj}\tag{SP}\\
 \mbox{s.t.}\quad&\epsilon_H ||\bm Y_{jt}||_2+\bm Y_{jt}^\mathsf{T}\hat{\bm{d}}^h-\hat S_{jt}^k\le M_h(1-\hat{w}^k_h),\forall j\in[J], t\in[T], h\in[T],\label{sp_cone2}\\
  \quad&Y_{ijt} \ge 0,~\forall i\in [I], j \in [J], t \in [T] \nonumber \\
  \quad& \eqref{con_demand},\eqref{con_u},\eqref{con_v},\eqref{con_z},\eqref{con_uvz},\eqref{dro_cone1}.\nonumber
\end{align}
\end{subequations}

If SP is feasible, the optimal objective value admits an upper bound of problem \eqref{dro_dflp_redef}; otherwise, we can obtain the optimal value of continuous decision variables according to the following SP$_\text{inf}$, which finds a solution that admits the lowest regret. As the constraint \eqref{dro_cone1} is always bounding at optimality, we mainly focus on the other second-order conic constraint \eqref{sp_cone2}. Introduce an auxiliary decision variable $r_{jth}\ge 0$, then,
\begin{subequations}
\begin{align}
 \quad \text{min.}\quad &\sum_{j\in[J]}\sum_{t\in[T]}\sum_{h\in[H]}r_{jth},\label{spinf_obj} \tag{SP$_\text{inf}$}\\
 \mbox{s.t.}\quad&\epsilon_H ||\bm Y_{jt}||_2+\bm Y_{jt}^\mathsf{T}\hat{\bm{d}}^h-\hat S_{jt}^k-r_{jth}\le M_h(1-\hat{w}^k_h),\forall j\in[J], t\in[T], h\in[T],\label{spinf_cone2}\\
 \quad&Y_{ijt} \ge 0,~\forall i\in [I], j \in [J], t \in [T] \nonumber \\
 \quad& \eqref{con_demand},\eqref{con_u},\eqref{con_v},\eqref{con_z},\eqref{con_uvz},\eqref{dro_cone1}.\nonumber
\end{align}
\end{subequations}

By solving SP or SP$_{\text{inf}}$ in the $k^{th}$ iteration, their optimal solutions ($\tilde{\bm Y}^k,\tilde {\bm u}^k$, $\tilde{\bm v}^k$, $\tilde{\bm z}^k$,$\tilde {\bm l}^k$) can help to generate valid OA cuts in MP.

\subsection{Master problem} \label{subsec:MP}
MP is an MILP with tailored OA cuts and provides a lower bound (LB) to the origin problem \eqref{dro_dflp_redef}. We first introduce the OA cuts added in the MP. Let $\Omega:=\{(\bm X,\bm S,\bm w)\in \{0,1\}^{J\times T}\times \mathbb Z^{J\times T}_+\times \{0,1\}^H: \eqref{SMLCPP},\eqref{dro_theta}\}$ denote the polyhedral set of integer variables, and $\Lambda:=\{(\bm Y,\bm u,\bm v,\bm z, \bm l)\in [0,1]^{I\times J\times T}\times \mathbb R^I\times \mathbb R^I\times \mathbb R^I \times \mathbb R^H:\{\eqref{con_demand},\eqref{con_u},\eqref{con_v},\eqref{con_z}\}\}$ denote the polyhedral set of continuous variables. Based on the convexity that has been proved in Lemmas \ref{lemma:convex1} and \ref{lemma:convex2}, any point $((\bm X,\bm S,\bm w),(\bm Y,\bm u,\bm v,\bm z, \bm l))\in \Omega \times \Lambda$, not necessarily feasible to  \eqref{dro_dflp_redef}, can generate the following OA
cuts for the \eqref{dro_dflp_redef}, see Propositions \ref{prop:OAcut1} and \ref{prop:OAcut2} for details.

\begin{proposition}\label{prop:OAcut1}
	Let $(\bar{\bm Y},\bar{\bm u},\bar{\bm v},\bar{\bm z},\bar{ \bm l})\in \Lambda$ be a feasible solution to the domain of continuous variables, then, the OA cut associated with constraint \eqref{dro_cone1} is
	\begin{equation}
		\label{OA_cut_cone1}
		\epsilon_H(\bar{\bm u}+\bar{\bm v}+\bar{\bm z})^\mathsf{T}(\bm u+\bm v+\bm z)+\left[(\bm u+\bm v+\bm z)^\mathsf{T}\bar{\bm{d}}^h -l_h\right]||\bar{\bm u}+\bar{\bm v}+\bar{\bm z}||_2 \le 0, \quad\forall h\in[H].
	\end{equation}
\end{proposition}

\begin{proof}
	For a given $h\in[H]$, the second-order conic inequality \eqref{dro_cone1} is equivalent to \[\phi_h(\bm u,\bm v,\bm z,l_h)=\epsilon_H||\bm u+\bm v+\bm z||_2+(\bm u+\bm v+\bm z)^\mathsf{T}\hat{\bm{d}}^h-l_h\le 0,\forall h\in[H].\] According to \eqref{def_phi}, by convexity of $\phi_h(\cdot)$, we have:
	\begin{equation}
		\phi_h(\bar{\bm u},\bar{\bm {v}},\bar{\bm {z}},\bar{l}_h)+\nabla \phi_h(\bar{\bm u},\bar{\bm v},\bar{\bm z},\bar{l}_h)^\mathsf{T}\left[ \bm u-\bar{\bm u},\bm v-\bar{\bm v},\bm z-\bar{\bm z},l_h-\bar{l}_h\right]^\mathsf{T}
		% \left[\begin{array}{c}
		% \bm u-\hat{\bm u}\\
		% \bm v-\hat{\bm v}\\
		% \bm z-\hat{\bm z}\\
		%  \eta_h-\hat{\bm \eta}_h
		% \end{array}\right]
		\le \phi_h({\bm u},{\bm v},{\bm z},{l_h})\le 0,
		\label{oa_cut_tmp1}
	\end{equation}
	where $\nabla \phi_h({\bm u},{\bm v},{\bm z},{l_h})=\left(\frac{\partial \nabla \phi_h({\bm u},{\bm v},{\bm z},{l_h})}{\partial\bm u},\frac{\partial \nabla \phi_h({\bm u},{\bm v},{\bm z},{l_h})}{\partial\bm v},\frac{\partial \nabla \phi_h({\bm u},{\bm v},{\bm z},{l_h})}{\partial\bm z},\frac{\partial \nabla \phi_h({\bm u},{\bm v},{\bm z},{l_h})}{\partial l_h}\right)$,
	$\frac{\partial \nabla \phi_h({\bm u},{\bm v},{\bm z},{l_h})}{\partial\bm u}=$$\hat{\bm{d}}^h + \\\frac{\epsilon_H\bm u}{||\bm u+\bm v+\bm z||_2}$, $\frac{\partial \nabla \phi_h({\bm u},{\bm v},{\bm z},{l_h})}{\partial\bm v}=\frac{\epsilon_H\bm v}{||\bm u+\bm v+\bm z||_2}+\hat{\bm{d}}^h$, $\frac{\partial \nabla \phi_h({\bm u},{\bm v},{\bm z},{l_h})}{\partial\bm z}=\frac{\epsilon_H\bm z}{||\bm u+\bm v+\bm z||_2}+\hat{\bm{d}}^h$ and $\frac{\partial \nabla \phi_h({\bm u},{\bm v},{\bm z},{l_h})}{\partial l_h}=-1$. As the Euclidean norm is no smaller than 0, inequality \eqref{oa_cut_tmp1} can be further simplified as \eqref{OA_cut_cone1} by plugging in the partial derivatives. If $\epsilon_H||\bar {\bm u}+\bar {\bm v}+\bar {\bm z}+\hat {l}_h||_2=0$, then, ignore the associated cut.
\end{proof}

\begin{proposition}\label{prop:OAcut2}
	Let $(\bar{\bm Y},\bar{\bm u},\bar{\bm v},\bar{\bm z},\bar{ \bm l})\in \Lambda$ be a feasible solution to the domain of continuous variables, then, the OA cut associated with constraint \eqref{dro_cone2} is,
	\begin{equation}
		\label{OA_cut_cone2}
		\epsilon_H\bar{\bm Y}_{jt}^\mathsf{T}\bm Y_{jt} +\left(\bm Y_{jt}^\mathsf{T}\bar{\bm{d}}^h-S_{jt}-M_h+M_hw_h\right)||\bar{\bm Y}_{jt}||_2\le 0, \quad\forall j\in [J], t\in [T], h\in[H].
	\end{equation}
\end{proposition}

\begin{proof}
	The proof is similar to that of Proposition \ref{prop:OAcut1}.
\end{proof}

%Obviously, in the $k^{th}$ iteration, the optimal solution $(\hat{\bm Y}^k,\hat {\bm u}^k$, $\hat{\bm v}^k$, $\hat{\bm z}^k$, $\hat {\bm l}^k)$ of MP and the optimal solution $(\tilde{\bm Y}^k,\tilde {\bm u}^k$, $\tilde{\bm v}^k$, $\tilde{\bm z}^k$, $\tilde {\bm l}^k)$ of SP both belong to set $\Lambda$. Then, by the convexity of continuous relaxation problem of the DR-MISOCP \eqref{dro_dflp_redef}, two types of cuts associated with these solutions can be added in MP:
Obviously, the optimal solutions of SP and SP$_\text{inf}$ belong to set $\Lambda$.  Then the OA cuts defined in Proposition \ref{prop:OAcut1} and Proposition \ref{prop:OAcut2} associated with these solutions can be added to the MP in the  $k^{th}$ iteration. Noting that the optimal solution to the ${(k-1)}^{th}$ iteration of MP, denoted as ($\hat {\bm Y}^{k-1}$, $\hat {\bm u}^{k-1}$, $\hat{\bm v}^{k-1}$, $\hat{\bm z}^{k-1}$, $\hat {\bm l}^{k-1}$),  also belongs to the set $\Lambda$. Consequently, we can also add the OA cuts associated this solution to MP in the $k^{th}$ iteration. 
\begin{subequations}
\begin{align}
\text{min.} \text{ } &\sum_{t\in[T]}\left\{f_{jt}X_{jt} + \sum_{j\in[J]} a_{jt}\left(S_{jt}-S_{j,t-1}\right)\right\}+\frac{1}{H}\sum_{h=1}^{H}l_h, \label{mp_obj}\tag{MP}\\
\mbox{s.t.} &\epsilon_H(\hat{\bm u}^\kappa\text{+}\hat{\bm v}^\kappa\text{+}\hat{\bm z}^\kappa)^\mathsf{T}(\bm u+\bm v+\bm z)+\left[(\bm u+\bm v+\bm z)^\mathsf{T}\bm \hat{\bm{d}}^h -l_h\right]||\hat{\bm u}^\kappa+\hat{\bm v}^\kappa+\hat{\bm z}^\kappa||_2 \le 0, \forall h\in[H],\kappa\in [k-1],\label{mp_cut1_1}\\
&\epsilon_H(\tilde{\bm u}^\kappa\text{+}\tilde{\bm v}^\kappa\text{+}\tilde{\bm z}^\kappa)^\mathsf{T}(\bm u+\bm v+\bm z)+\left[(\bm u+\bm v+\bm z)^\mathsf{T}\hat{\bm{d}}^h -l_h\right]||\tilde{\bm u}^\kappa+\tilde{\bm v}^\kappa+\tilde{\bm z}^\kappa||_2 \le 0, \forall h\in[H],\kappa\in [k],\label{mp_cut1_2}\\
& \epsilon_H\bm Y_{jt}^\mathsf{T}\hat{\bm Y}_{jt}^\kappa +\left(\bm Y_{jt}^\mathsf{T}\hat{\bm{d}}^h-S_{jt}-M_h +M_hw_h\right)||\hat{\bm Y}^\kappa_{jt}||_2\le 0, \forall j\in [J], t\in [T], h\in[H],\kappa\in [k-1],\label{mp_cut2_1}\\
& \epsilon_H\bm Y_{jt}^\mathsf{T}\tilde{\bm Y}_{jt}^\kappa +\left(\bm Y_{jt}^\mathsf{T}\hat{\bm{d}}^h-S_{jt}-M_h+M_hw_h\right)||\tilde{\bm Y}^\kappa_{jt}||_2\le 0, \forall j\in [J], t\in [T], h\in[H],\kappa\in [k],\label{mp_cut2_2}\\
\text{ } & \eqref{Trival_Con},\eqref{con_u}\text{-}\eqref{con_uvz},\eqref{dro_theta},\eqref{dro_theta_domain}.\nonumber
\end{align}
\end{subequations}

\begin{remark}\label{my_remark_oa}
	If $k = 0$, then we only add the OA cuts associated with the optimal solution to the SP problem, as we have not solved an MP problem yet.
\end{remark}

\subsection{Procedure}\label{subsec:procedure}
According to \cite{1987An}, OA algorithm operates by iteratively solving SP and MP, updating UB and LB respectively. On the basis of the formulations of MP and SP, the procedure of the proposed OA algorithm can be summarized as follows.

\vspace{0.1in}
\begin{algorithm}[H]
	\SetAlgoLined
	{\bf Input}: ($\hat S_{jt}^0, \hat X_{jt}^0, \hat w_h^0$): the initial values of integer variable; K: maximum number of iterations, indexed by $k$, initiate $k = 0$\;
	Solve SP with fixed integer values ($\hat S_{jt}^0, \hat X_{jt}^0, \hat w_h^0$), obtain the incumbent solution $(\tilde{\bm Y}^k, \tilde{\bm u}^k,\tilde{\bm v}^k,\tilde{\bm z}^k,\tilde{\bm l}^k)$. let UB$^0$ be the optimal value of SP\;
	Build model MP, add objective, linear constraints \eqref{Trival_Con}, \eqref{con_u}\text{-}\eqref{con_uvz}, \eqref{dro_theta}, \eqref{dro_theta_domain}, OA cuts associated with $(\tilde{\bm Y}^k, \tilde{\bm u}^k,\tilde{\bm v}^k,\tilde{\bm z}^k,\tilde{\bm l}^k)$ to initialize MP\;
	Solve MP, get the optimal value of continuous variables $(\hat{\bm Y}^k, \hat{\bm u}^k,\hat{\bm v}^k,\hat{\bm z}^k,\hat{\bm l}^k)$ and integer variables ($\hat S_{jt}^k, \hat X_{jt}^k, \hat w_h^k$), and let LB$^k$ be the optimal objective value of MP\;
	k = k + 1 \;
	\For{$ k = 1:K$}{
		Solve SP with fixed integer variables ($\hat S_{jt}^k, \hat X_{jt}^k, \hat w_h^k$) \;
		\uIf{SP is feasible}
		{
			Obtain an incumbent solution $(\tilde{\bm Y}^k, \tilde{\bm u}^k,\tilde{\bm v}^k,\tilde{\bm z}^k,\tilde{\bm l}^k)$; let UB$^{k}$ be the optimal objective of SP\;
		}
		\uElse
		{Solve SP$_{\text{inf}}$, obtain the incumbent solution $(\tilde{\bm Y}^k, \tilde{\bm u}^k,\tilde{\bm v}^k,\tilde{\bm z}^k,\tilde{\bm l}^k)$; let UB$^{k}$=UB$^0$ \;}
		Add OA cuts associated with $(\hat{\bm Y}^{k-1}, \hat{\bm u}^{k-1},\hat{\bm v}^{k-1},\hat{\bm z}^{k-1},\hat{\bm l}^{k-1})$ and $(\tilde{\bm Y}^{k-1}, \tilde{\bm u}^k,\tilde{\bm v}^k,\tilde{\bm z}^k,\tilde{\bm l}^k)$ to MP \;
		Solve MP, get the optimal value of continuous variables $(\hat{\bm Y}^k, \hat{\bm u}^k,\hat{\bm v}^k,\hat{\bm z}^k,\hat{\bm l}^k)$ and integer variables ($\hat S_{jt}^k, \hat X_{jt}^k, \hat w_h^k$), and let LB$^k$ be the optimal objective value of MP\;
		\uIf{The gap between current LB and UB is smaller than a threshold}
		{break;}
	}
	{\bf Output}: optimal objective values UB, LB, the optimal location and capacity ($\hat {\bm X}^k$,$\hat {\bm S}^k$)\;
	\caption{An outer approximation algorithm of the problem \eqref{dro_dflp_redef}} \label{algo_OA}
\end{algorithm}

\section{Numerical results}\label{section:numerical}
We conduct numerical experiments to validate the performance of our MFLCP-W model and the proposed OA algorithm in this section. All experiments are implemented on a Ubuntu operating system with AMD Ryzen Threadripper 3990X 64-Core, 128-Thread Unlocked Desktop Processor. The CPLEX 12.8 optimization toolbox for MATLAB is employed to solve the MP and SP in the OA algorithm.

\subsection{Performance of the proposed MFLCP-W model }\label{subsection:performance}
Numerical experiments on a randomly generated $10\times10$ complete graph, representing the affected areas and facility candidates, are conducted to reveal the performance of the MFLCP-W model \eqref{dro-p}. Parameters similar to those in \citet{liu2019distributionally,zhang2015novel} are summarized in Table \ref{para_set}.
% \begin{longtable}{ll}
% \label{Parameter}
% \endfirsthead
% \endhead
%  $T$     & Equals to 3. \\
%  $f_{jt}$& Uniformly generated from $[100(T-t),100+100(T-t)]$ in period $t$.\\
%  $a_{jt}$& Uniformly generated from $[0,2]$ in period $t$. \\
%  $c_{ji}$& Uniformly generated from $[0,5]$.\\
%  $q_t$ &Equals to $20t$ in period $t$.\\
%  $d_i$ &Uniformly generated from $[0,30]$.\\
%  $\alpha$ &Equals to 0.8.\\
%  $\eta$ &Equals to 0.8.
% \end{longtable}
\renewcommand{\arraystretch}{1} %¿ØÖÆ±í¸ñÐÐ¸ßµÄËõ·Å±ÈÀý
\begin{table}[!h]
	\caption{Parameters Setting}
	%	\fontsize{8}{10}\selectfont
	\begin{tabu} to 1\textwidth {X[1,l] X[12,l] }
		\toprule[1 pt]
		Parameter & Setting  \\ \midrule
		$T$     & Equals to 3. \\
        $f_{jt}$& Uniformly generated from $[100(T-t),100+100(T-t)]$ in period $t$.\\
        $a_{jt}$& Uniformly generated from $[0,2]$ in period $t$. \\
        $c_{ji}$& Uniformly generated from $[0,5]$.\\
        $q_t$ &Equals to $20t$ in period $t$.\\
        $d_i$ &Uniformly generated from $[0,30]$.\\
        $\alpha$ &Equals to 0.8.\\
         $\eta$ &Equals to 0.8. \\
		\bottomrule[1 pt]		
	\end{tabu}
	\label{para_set}	
\end{table}

As it is difficult and expensive to install large-scaled facilities right after an emergency, we assume the initial construction cost is significantly larger than that in other periods, and the available construction capacity also increases steadily along with time.
\noindent Penalty cost is generated by the definition of DLF for medicines, i.e., $\mu_t=\frac{9.772697}{1+3.9031\cdot e^{-0.7919 t}}$ \citep{Wang2017poms}. %Sensitivity analysis associated with different forms of penalty cost is conducted in Subsection \ref{subsection:sensitivity}.

We evaluate the impacts of Wasserstein radius $\epsilon_H$ on the out-of-sample performance of our model, and compare the proposed model with the SAA methods in terms of the out-of-sample performance, which is quantified by examining the objective values \eqref{out_of_sample_obj} and  satisfaction probability  \eqref{out_of_sample_con} under \textit{new} samples, i.e.,
\begin{equation}
\label{out_of_sample_obj}
  \sum_{t\in[T]}\left\{f_{jt}X_{jt} + \sum_{j\in[J]} a_{jt}\left(S_{jt}-S_{j,t-1}\right)\right\} +  \mathbb{E}_{\mathbb{F}}\left[\bm u^\mathsf{T}\bm{d} +\bm v^\mathsf{T}\bm{d}+\bm{z}^\mathsf{T}\bm{d} \right],
\end{equation}

\begin{equation}
\label{out_of_sample_con}
  \mathbb{P}_{\mathbb{F}}\left\{\max_{j \in [J], t \in [T]}\{\bm{Y}_{jt}^\mathsf{T}\bm{d}-S_{jt}\} \le 0 \right\}.
\end{equation}

It is difficult to calculate the objective and probability due to the multivariate integral, we instead generate $H_T=100$ samples from the uniform distribution as testing samples to approximate \eqref{out_of_sample_obj} and \eqref{out_of_sample_con}, i.e.,

\begin{align}
 \label{OBJ}
  \sum_{t\in[T]}\left\{f_{jt}X_{jt} + \sum_{j\in[J]} a_{jt}\left(S_{jt}-S_{j,t-1}\right)\right\} + \frac{1}{H_T}\sum_{h=1}^{H_T}\left[\bm u^\mathsf{T}\hat{\bm{d}}^h +\bm v^\mathsf{T}\hat{\bm{d}}^h+\bm{z}^\mathsf{T}\hat{\bm{d}}^h \right], \tag{objective}
\end{align}
\begin{align}
 \label{probability}
  \frac{1}{H_T}\sum_{h=1}^{H_T}\left\{\max_{j \in [J],t \in [T]}\left\{ \bm{Y}_{jt}\bm{\hat{d}}^h - S_{jt} \right\}\le 0\right\}, \tag{probability}
\end{align}
where $\hat{\bm d}^h$ is the $h^\text{th}$ testing sample.

\subsubsection{Impact of the Wasserstein ball radius}\label{subsection:radius}
\begin{figure}[ht] \centering
\subfigure[H = 20] { \label{fig:a}
\includegraphics[width=0.3\linewidth]{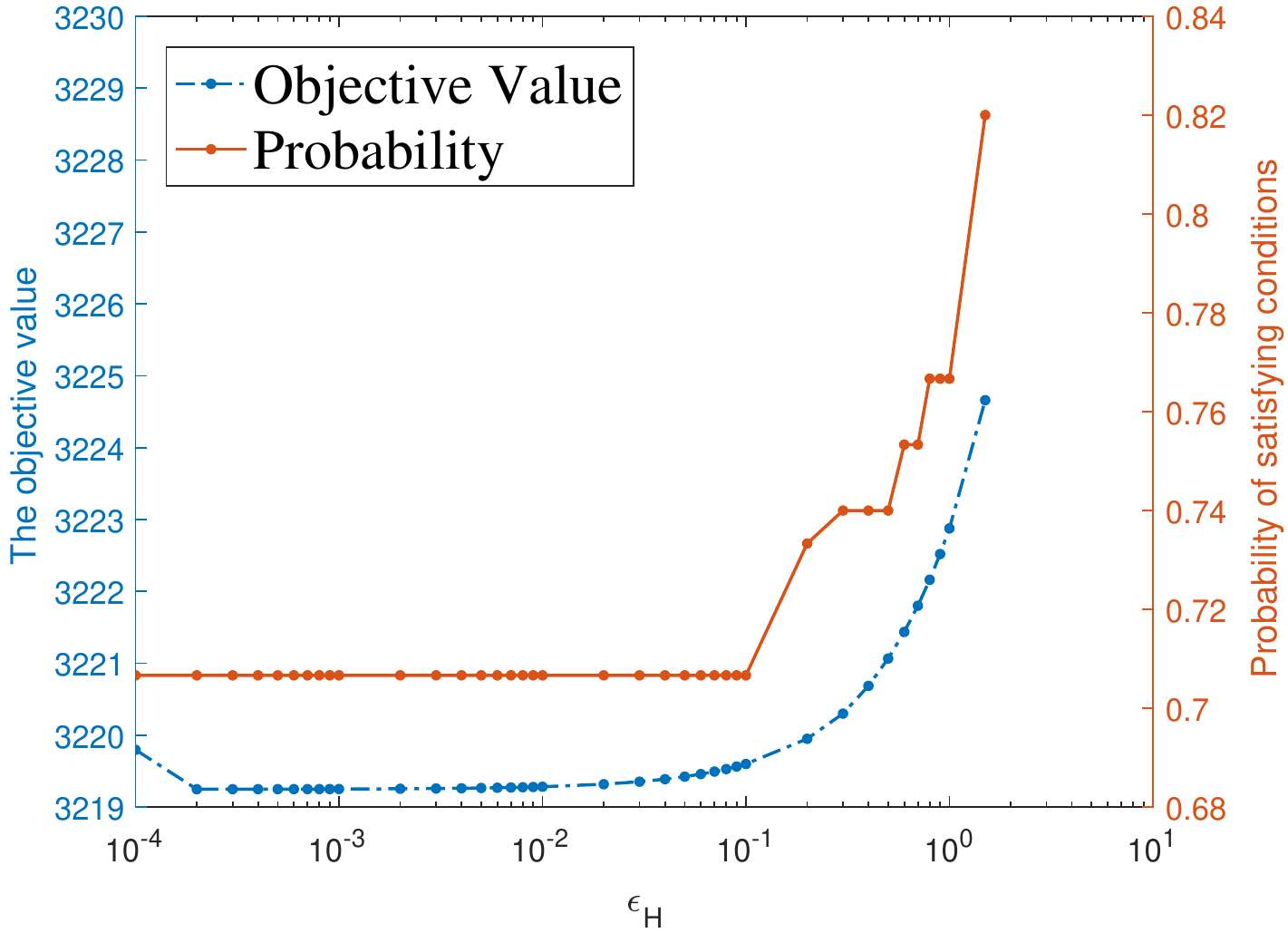}%{cvN=20_k=30.pdf}
}
\subfigure[H = 50] { \label{fig:b}
\includegraphics[width=0.3\linewidth]{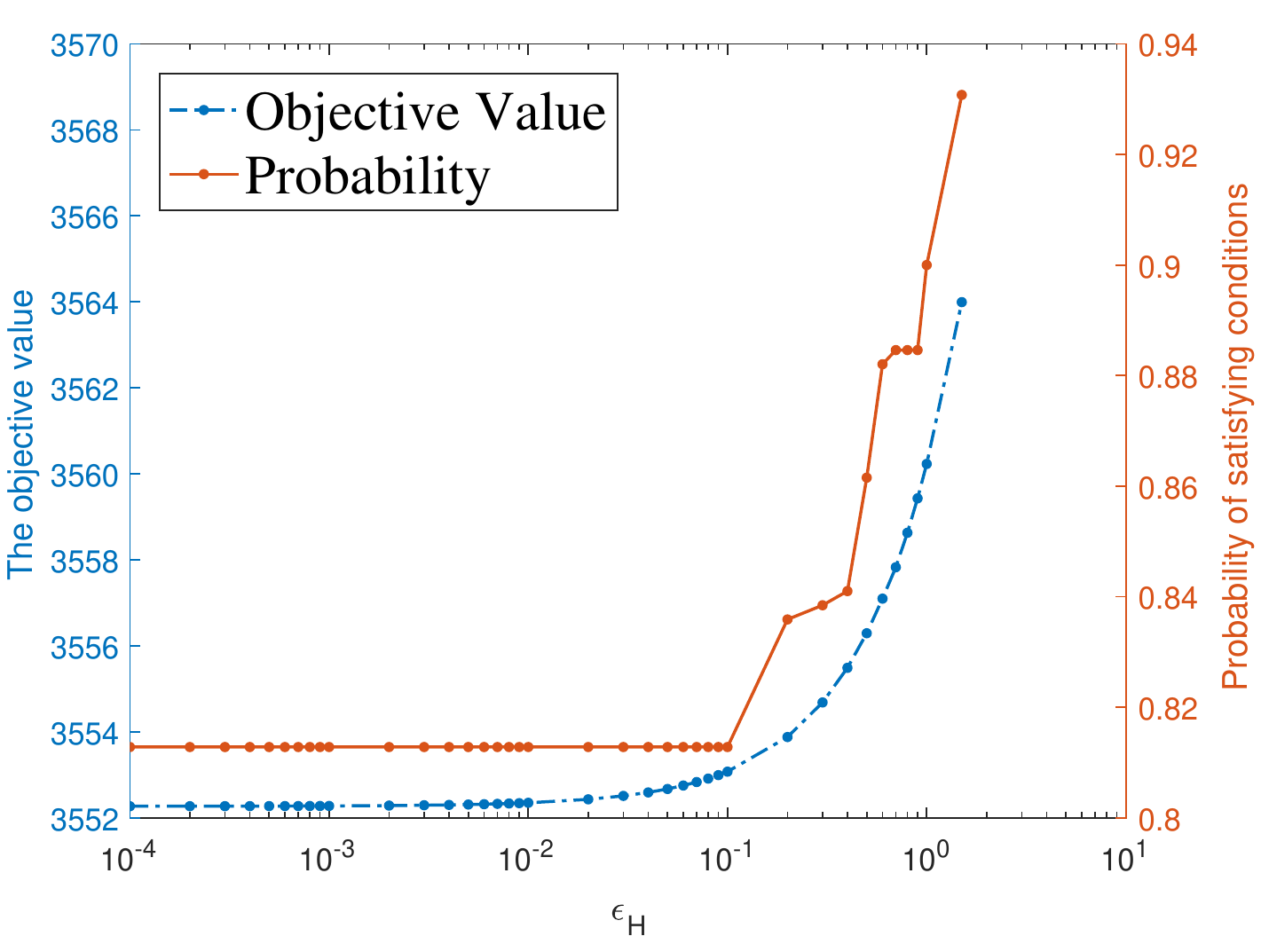}%{cvN=50_k=30.pdf}
}
\subfigure[H = 80] { \label{fig:c}
\includegraphics[width=0.3\linewidth]{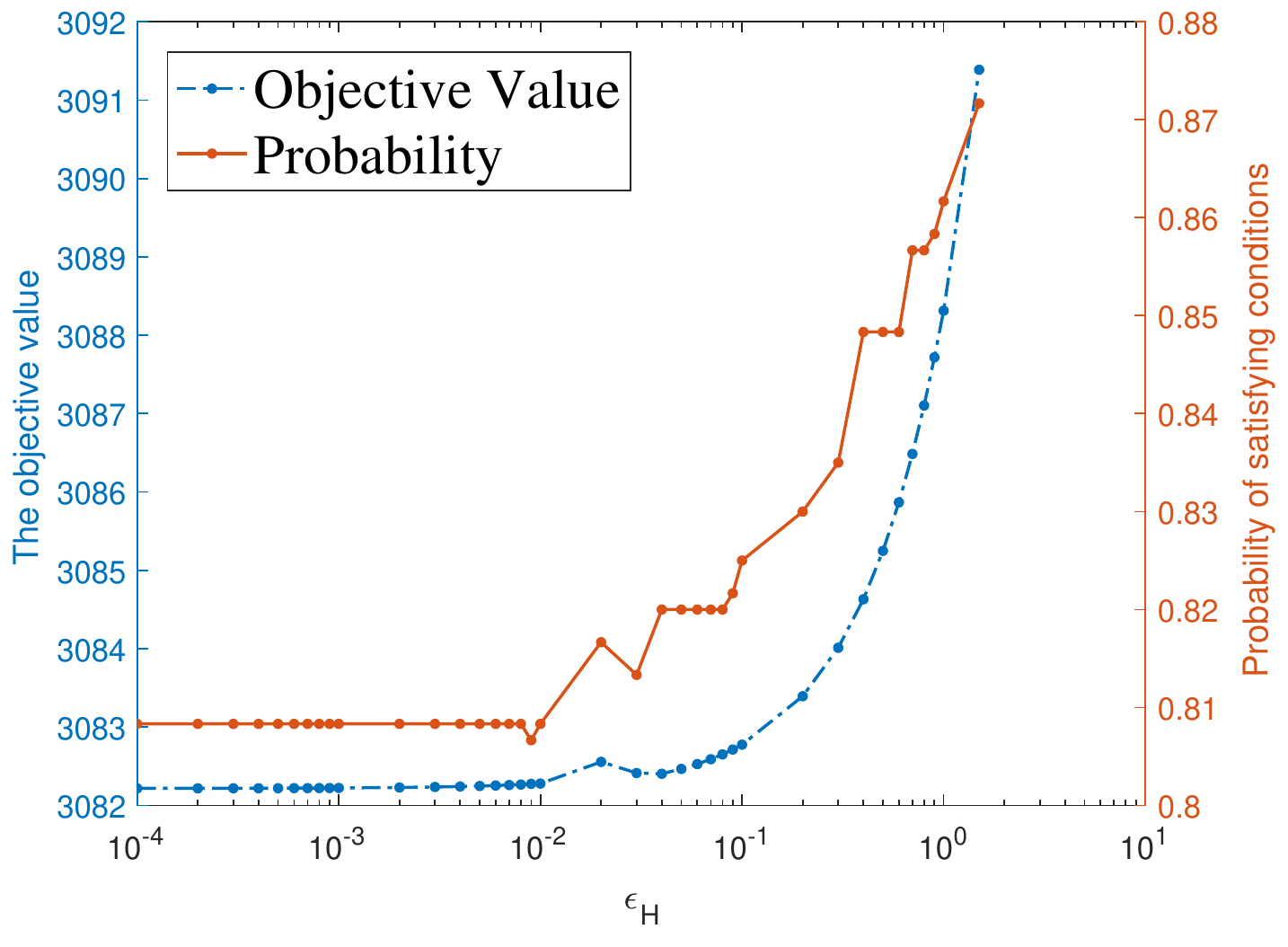}%{cvN=80_k=30.pdf}
}
\caption{ The out-of-sample obj (left axis, blue line) and probability (right axis, orange line) as a function of the Wasserstein radius $\epsilon_H$ and averaged over $30$ experiments. (a) H = 20. (b)H = 50. (c) H = 80.}
\label{fig:radius_valida}
\end{figure}

To evaluate the impact of radius $\epsilon_H$ on the proposed MFLCP-W model, we solve problem \eqref{dro-equi-socp} using training datasets with the cardinality $H \in \{20,50,80\}$. %The Wasserstein radii $\epsilon_H$ are drawn from set $\mathcal{E}$.% $\{0.0001,0.0002,\dots,1,1.5\}$
The averaged out-of-sample \ref{OBJ} and \ref{probability} under 30 independent experiments are shown in Figure \ref{fig:radius_valida}. From Figure \ref{fig:a}, the out-of-sample \ref{OBJ} firstly improves (decreases) up to a critical Wasserstein radius, then deteriorates (increases) as the radius rises when $H$ is relatively small; and \ref{probability} is non-decreasing in $\epsilon_H$. This observation has a intuitive appeal since the probability of including true distribution $\mathbb{F}\in \mathcal{F}^\infty_H$ generally increases as $\epsilon_H$ grows, and hence joint chance constraints can be satisfied more easily. However, an excessively large radius $\epsilon_H$ increases the conservatism of the proposed model, and leads to a worse out-of-sample objective. Hence, we need select a proper radius to obtain solutions that balance the values of out-of-sample satisfaction probability and the objective value, see Remark \ref{my_remark} for more details.

\subsubsection{Comparison with the state-of-the-art methods}
\label{subsection:saa}
This subsection compares the proposed model with SAA in terms of the out-of-sample \ref{OBJ} and \ref{probability}. We set the confidence level $\eta = 0.2$ and the sample dataset size $H = \{10,20,\dots,100\}$. Ideally, one should select the radius which achieves the optimal out-of-sample performance over all $\epsilon_H > 0$. However, it is impossible as the true distribution $\mathbb{F}$ is unknown. Instead, we consider a discrete set $\mathcal{E}$ with cardinality $E$,
$$\mathcal{E} = \left\{\epsilon^1_H,\dots,\epsilon^E_H\right\} =\{0.0001:0.0001:0.001, 0.002:0.001:0.01,0.02:0.01:0.1,0.2:0.1:1,1.5\},$$
and select the optimal radius by the following holdout method:
\begin{enumerate}[.]
 \item[$\bullet$] \textit{Holdout method:} Randomly select $80\%$ samples from a dataset $\{\hat{\bm d}^1,\dots,\hat{\bm d}^H\}$ as the training set and the remaining $20\%$ as the validation set. First, use the training set to solve problem \eqref{dro_dflp_redef} under all candidate radii and obtain the corresponding optimal solutions. Then, use the validation set to estimate the out-of-sample \ref{OBJ} and \ref{probability} of the optimal solution. Finally, the optimal radius $\epsilon_H$ is set as the one with the best out-of-sample performance based on Remark \ref{my_remark}.
\end{enumerate}

\begin{remark}\label{my_remark}
Since it is impossible to find a radius simultaneously exhibiting the best out-of-sample \ref{OBJ} and \ref{probability}, we seek a trade-off between the two objectives. Let $OBJ(\epsilon^i_H)$ and $Pro(\epsilon^i_H)$ be the out-of-sample objective value and satisfaction probability,  $\forall \epsilon^i_H \in \mathcal{E}, i \in [E]$. The optimal radius $\epsilon_H$ in the holdout method is selected as,
\begin{equation*}
 \epsilon_H \in \arg\max_{\epsilon^i_H \in \mathcal{E}} \left\{ 1 - \frac{OBJ(\epsilon^i_H)}{\sum_{i=1}^{E}OBJ(\epsilon^i_H)} + \frac{Pro(\epsilon^i_H)}{\sum_{i=1}^{E}Pro(\epsilon^i_H)} \right\}.
\end{equation*}
\end{remark}

\begin{figure}[ht] \centering
\subfigure { \label{fig:saa_con}
\includegraphics[width=0.45\linewidth]{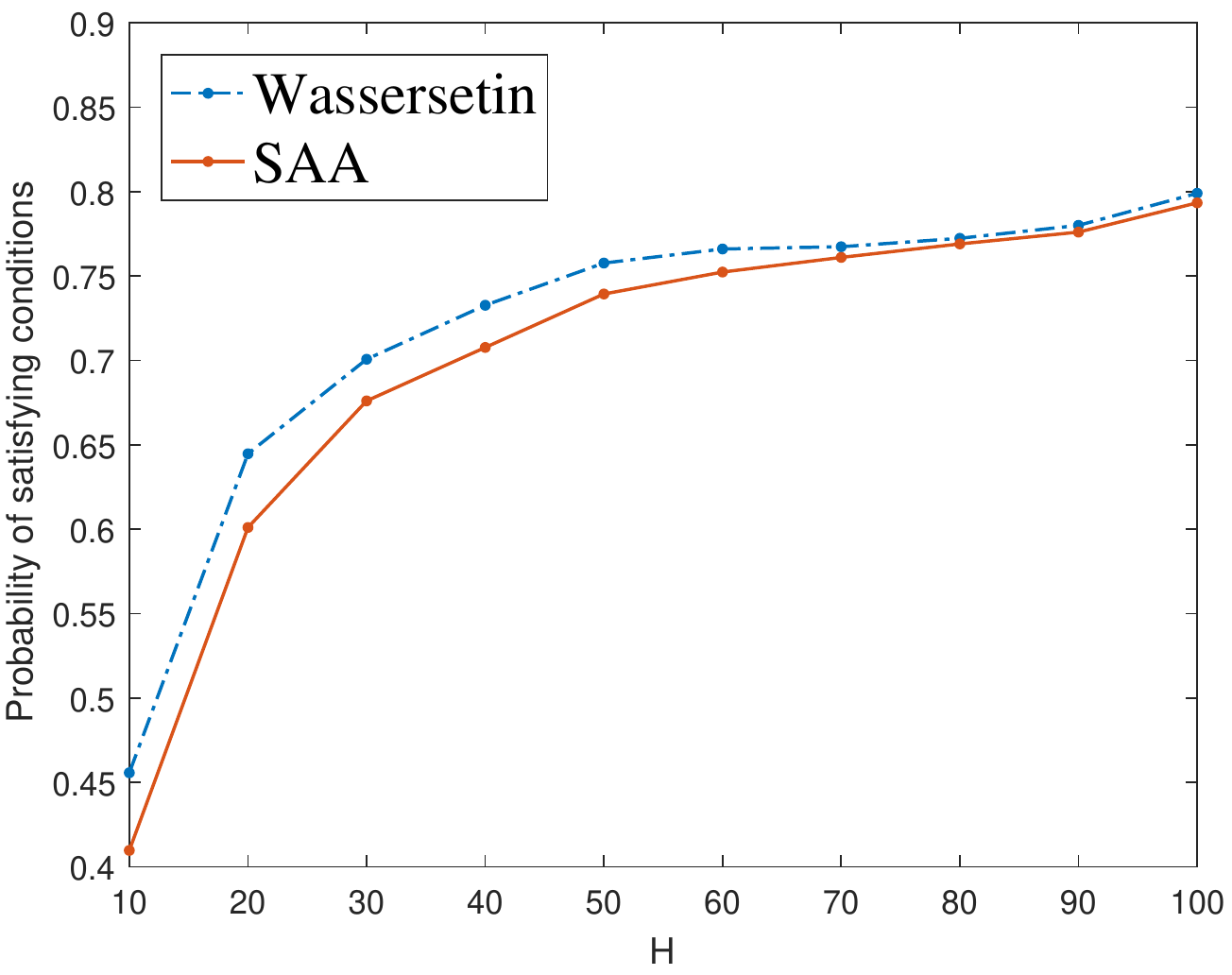}
}
\subfigure { \label{fig:saa_obj}
\includegraphics[width=0.45\linewidth]{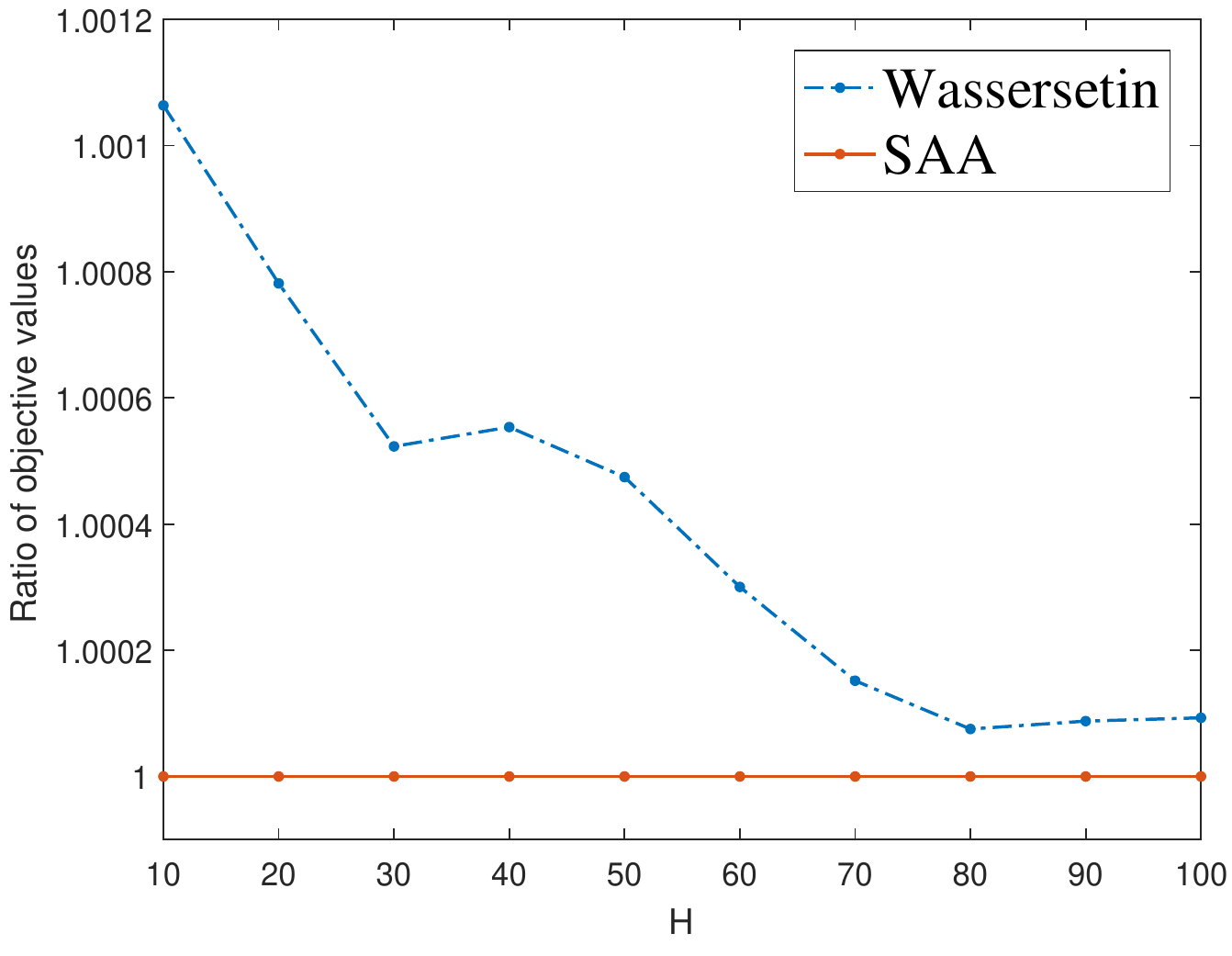}
}
\caption{ The left graph shows the out-of-sample \ref{probability} as a function of the sample dataset size $H$ of the SAA model (red) and the MFLCP-W model (blue). The right graph compares of the out-of-sample \ref{OBJ} of the two models.}
\label{fig_saa}
\end{figure}

We solve the proposed MFLCP-W model and the SAA model under different sample sizes $H$, and evaluate their out-of-sample \ref{OBJ} and \ref{probability} respectively. Results in Figure \ref{fig_saa} imply that the proposed MFLCP-W model produces solutions with better out-of-sample satisfaction probability, at a very modest objective value increases (typically less than $0.1\%$). In HL, system reliability and relief effectiveness are more important than the cost and consequently the MFLCP-W is more suitable to the relief process. Moreover, Figure \ref{fig_saa} shows that the optimal radius follows a decreasing trend with the increment of sample size, which is consistent with the convergence trend of the Wasserstein ambiguity set, i.e., the Wasserstein ball can asymptotically degenerate to the true distribution as the sample size increases to infinity \citep{esfahani2018data}. 

\subsection{Performance analysis of the OA algorithm}\label{subsection:cplex}
The goal of this subsection is to make a comprehensive comparison of the OA algorithm and CPLEX. Twenty combination of scale parameters ($I=J$, $T$ and $H$) with five instances for each are conducted. The algorithms are terminated once total computational time is larger than 1800 seconds.%, or the MIP gap reaches 0.001.

Table \ref{tab:cplex} reports the numerical results as three aspects: total computational time (Time), mixed integer gap within the time limit (Gap), and the percentage of getting the optimal solution within 1800 seconds (Successfully solved ratio), where the first two metrics are reflected by the average value of five instances. Columns ``C" and ``OA" represent the performance of CPLEX and the proposed OA method respectively.

\begin{table}[!htb]
 \centering
 \caption{Performance analysis of the OA algorithm}
 \begin{tabular}{ccccccccc}
 \toprule
 \multirow{2}[0]{*}{$I$=$J$} & \multirow{2}[0]{*}{$T$} & \multirow{2}[0]{*}{$H$} & \multicolumn{2}{c}{Time} & \multicolumn{2}{c}{Gap} & \multicolumn{2}{c}{Successfully solved ratio} \\
\cmidrule(r){4-5}\cmidrule(r){6-7}\cmidrule(r){8-9}
   &  &  & C  & OA & C  & OA & C  & OA \\ \midrule
 5  & 3  & 20 & 0.39 & 0.36 & 0  & 0  & 100\% & 100\% \\
 5  & 3  & 50 & 1.79 & 1.13 & 0  & 0  & 100\% & 100\% \\
 5  & 3  & 100 & 31.20 & 28.90 & 0  & 0  & 100\% & 100\% \\
 15 & 3  & 20 & 11.55 & 6.75 & 0  & 0  & 100\% & 100\% \\
 15 & 3  & 50 & 194.91 & 121.32 & 0  & 0  & 100\% & 100\% \\
 15 & 3  & 100 & 1800 & 1800 & 3.56\% & 0.17\% & 0\% & 0\% \\
 30 & 3  & 20 & 45.73 & 43.32 & 0  & 0  & 100\% & 100\% \\
 30 & 3  & 50 & 1627.98 & 1036.39 & 1.67\% & 0.03\% & \textbf{20\%} & \textbf{80\%} \\
 30 & 3  & 100 & 1800 & 1800 & 6.56\% & 0.13\% & 0\% & 0\% \\
 40 & 3  & 20 & 117.73 & 98.99 & 0  & 0  & 100\% & 100\% \\
 5  & 5  & 20 & 0.81 & 0.69 & 0  & 0  & 100\% & 100\% \\
 5  & 5  & 50 & 8.16 & 6.59 & 0  & 0  & 100\% & 100\% \\
 5  & 5  & 100 & 32.15 & 30.86 & 0  & 0  & 100\% & 100\% \\
 15 & 5  & 20 & 32.78 & 27.31 & 0  & 0  & 100\% & 100\% \\
 15 & 5  & 50 & 1800 & 1034.87 & 1.43\% & 0.07\% & \textbf{0\%} & \textbf{60\%} \\
 15 & 5  & 100 & 1800 & 1800 & 6.04\% & 0.17\% & 0\% & 0\% \\
 30 & 5  & 20 & 259.33 & 130.94 & 0  & 0  & 100\% & 100\% \\
 30 & 5  & 50 & 1800 & 1800 & 3.97\% & 0.13\% & 0\% & 0\% \\
 30 & 5  & 100 & 1800 & 1800 & 10.47\% & 0.13\% & 0\% & 0\% \\
 40 & 5  & 20 & 432.07 & 309.26 & 0  & 0  & 100\% & 100\% \\
  \bottomrule
 \end{tabular}%
 \label{tab:cplex}%
\end{table}%

From Table \ref{tab:cplex}, the average computational time of OA is generally smaller than that of CPLEX, especially in large-scaled cases, e.g., $I=J=30$, $T=5$ and $H=20$, where the computational time of OA takes only 50\% of CPLEX. For those instances that cannot be solved within the time limit, OA achieves a smaller optimality gap. For example, when $I=J=30$, $T=5$ and $H=100$, although both algorithms fail to get any optimal solution, the gap of OA (0.13\%) is much smaller than that of CPLEX (10.47\%); when $I=J=30$, $T=3$ and $H=50$, CPLEX only computes one instance to optimality, while OA gets four optimal solutions out of five.

\subsection{Sensitivity analysis}\label{subsection:sensitivity}
In this section, we test the impact of different forms of unit penalty cost $\mu(t)$ in our MFLCP-W model. Details of the selected functions are summarized as follows. The first two are the mixed logit models utilized in \cite{cantillo2018discrete}, where one is specified with Box-Cox transformations $\mu_B(t)$, and the other is based on the exponential transformation $\mu_E(t)$. The third function $\mu_L(t)$ is an S-shape logistic growth function adopted in \cite{Wang2017poms}, which naturally converges to a specific maximum value $\text{L}_{\max}$. To make a fair comparison, we normalize the values of $\mu_B(t)$ and $\mu_E(t)$ to the range of $\mu_L(t)$, so that the transformed functions lie in the same order of magnitude. Moreover, as the DLF is specified in days, we transform the unit of deprivation time to hours by dividing 24. 

\begin{figure}[htbp] \centering
\subfigure[Illustration of penalty cost functions.]{ \label{fig:fig_penalty}
\includegraphics[width=0.48\linewidth]{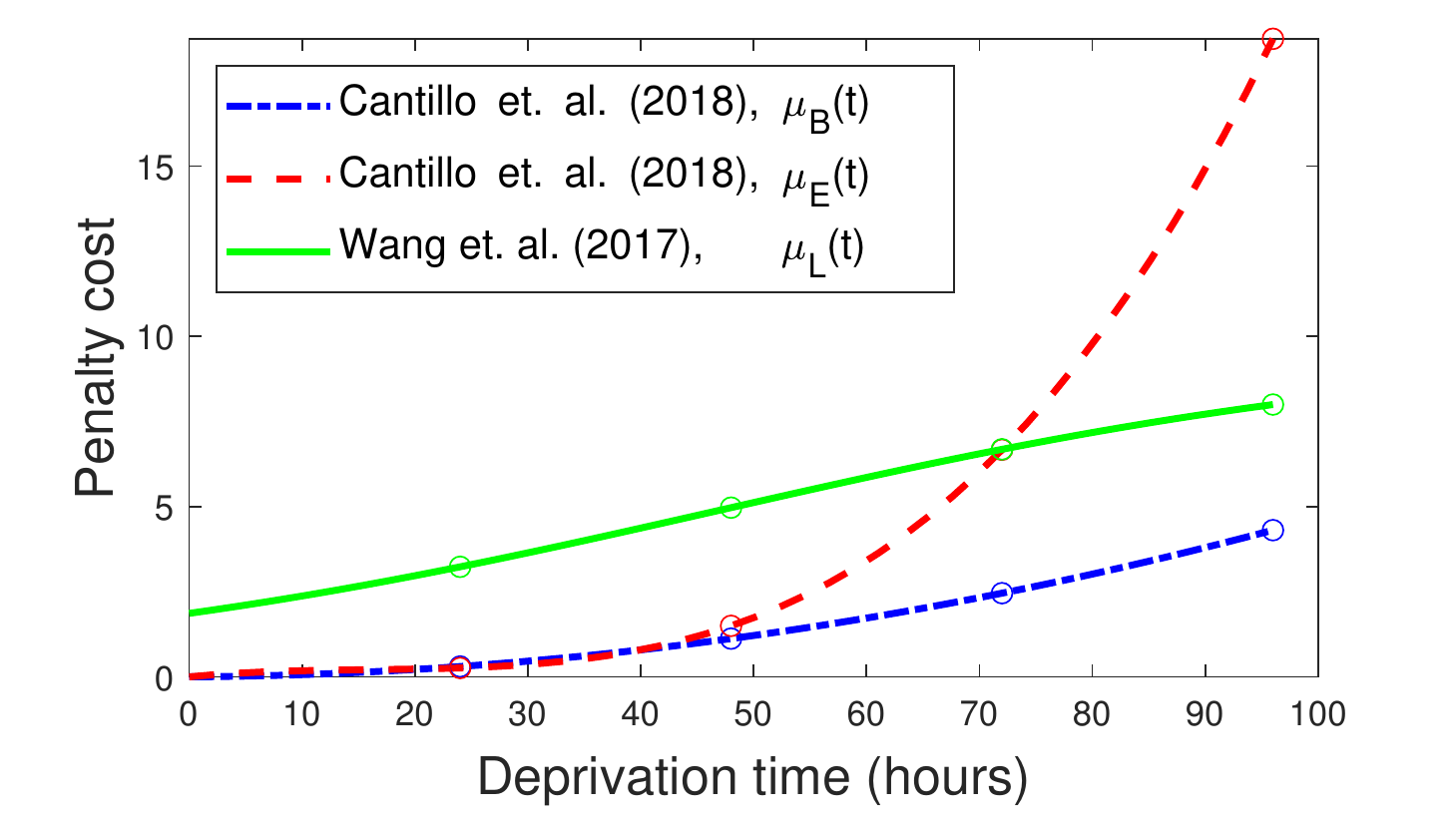} %[]
}
\subfigure[Values of penalty costs at the end of the day.]{ \label{fig:tab_penalty}
\includegraphics[width=0.45\linewidth]{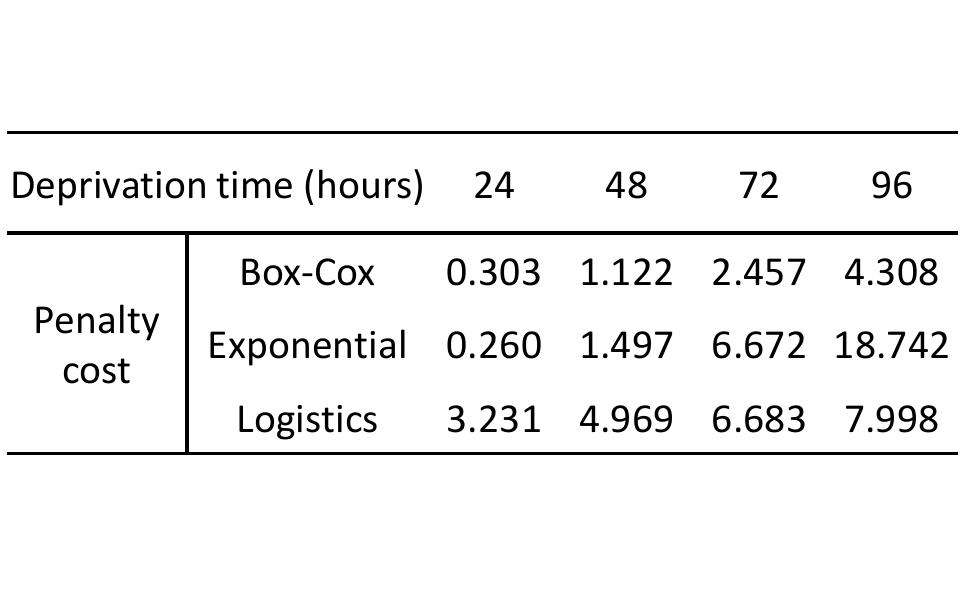} %
}
\caption{Unit penalty cost}
\label{fig_penalty}
\end{figure}

In summary, three types of penalty cost,  are illustrated in equations \eqref{mub}-\eqref{mul} and Figure \ref{fig:fig_penalty}.
\begin{subequations}
\begin{align}
    \mu_B(t) &=(459.52t^2+1935.5)t/2000\cdot\frac{\text{L}_{\max}}{\text{M}_{\max}},\label{mub}\\
    \mu_E(t) &=(36.606t^3-1764.7t^2+32375t)/2000\cdot\frac{\text{L}_{\max}}{\text{M}_{\max}},\label{mue}\\
    \mu_L(t) &=   \frac{9.745492}{1+4.228e^{-0.7407\cdot (t/24)}}, \label{mul}
\end{align}
\end{subequations}
where $\text{M}_{\max}$ is a given constant. By definition of \eqref{mul}, it is obviously that $\text{DLF}_{\max}$ equals to $9.745492$.

Figure \ref{fig:fig_penalty} shows that the values of penalty cost gradually increase with the planning horizon. The $\mu_E(t)$ achieves the largest growth rate, and the other two slopes gently. Note that, $\mu_L(t)$ takes much larger value than the others during the initial hours right after an emergency. In the experiments, let the planning horizon $T=96$ hours, unit penalty cost $\mu(t)$ equal to the values of Table \ref{fig:tab_penalty}, the scale parameters be $I=J=30$, $H=20$ and other parameters take the same value of those in Section \ref{subsection:performance}. The graph is also randomly generated as in \ref{subsection:performance}.

\begin{figure}[htbp]
    \centering
    \includegraphics[width=1.0\linewidth]{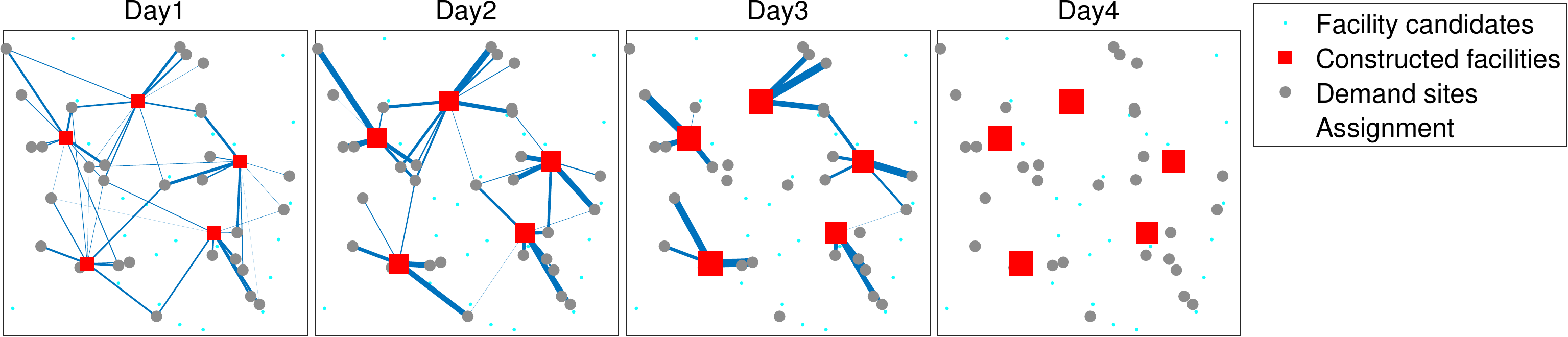}
    \caption{Topology analysis of the model with unit penalty cost $\mu_L(t)$.}
    \label{fig:logistics}
\end{figure}

Table \ref{tab:sensitivity} reports the numerical results of different penalty cost. From Table \ref{tab:sensitivity}, we observe that the $\mu_E(t)$ tends to open more facilities, which incurs a larger total logistics cost and a smaller penalty cost; $\mu_L(t)$ chooses to open smaller number of facilities and expand their capacities in the later periods (see Figure \ref{fig:logistics}), however, the strategy may cause large amount of penalty cost during the entire planning horizon.

\begin{table}[htbp]
  \centering
  \caption{Comparison of different types of penalty cost}
    \begin{tabular}{lllcc}\toprule
    \multicolumn{2}{c}{Penalty cost} & The index of established nodes & Total logistics cost & Total penalty cost \\\midrule
    Box-box & $\mu_B(t)$  &  2, 3, 5, 7, 10, 16 & 1419.09 & 1152.38 \\
    Exponential& $\mu_E(t)$&  2, 3, 5, 7, 10, 13, 16, 22, 25 & 1933.31 & 929.42 \\
    Logistics & $\mu_L(t)$& 2, 3, 7, 10, 16 & 1250.42 & 2323.47 \\\bottomrule
    \end{tabular}%
  \label{tab:sensitivity}%
\end{table}%

\begin{figure}[htbp]
    \centering
    \includegraphics[width=1.0\linewidth]{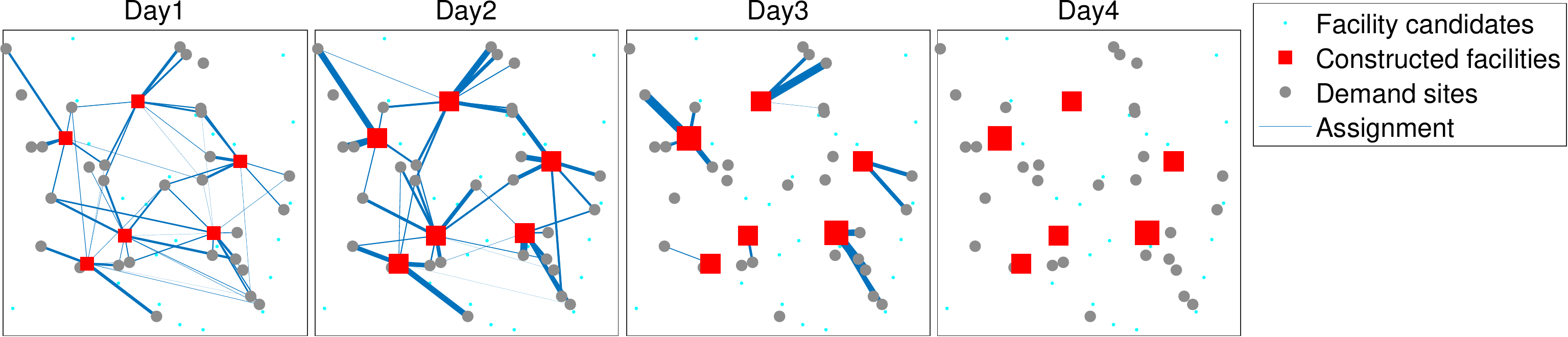}
    \caption{Topology analysis of the model with unit penalty cost function $\mu_B(t)$.}
    \label{fig:box}
\end{figure}

\begin{figure}[htbp]
    \centering
    \includegraphics[width=1.0\linewidth]{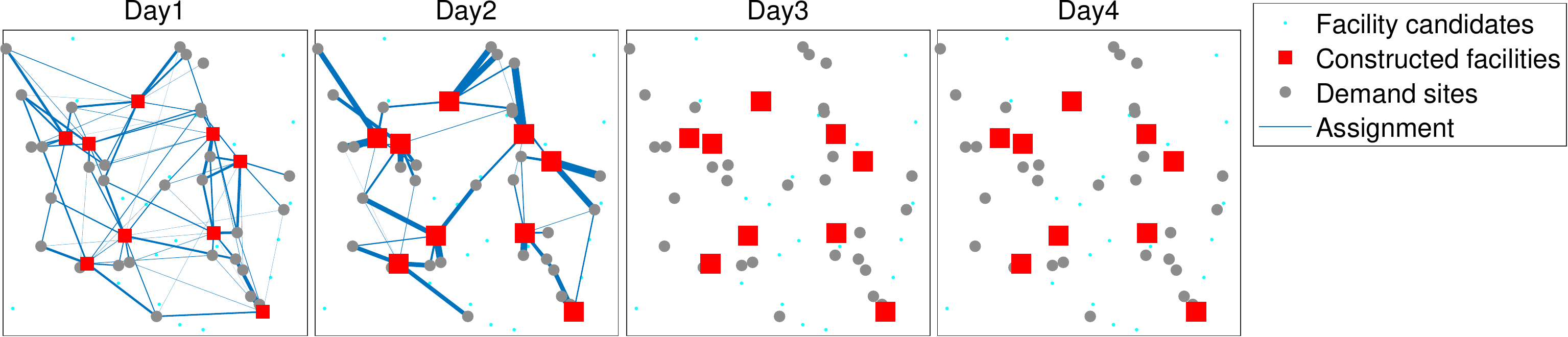}
    \caption{Topology analysis of the model with unit penalty cost function $\mu_E(t)$.}
    \label{fig:exp}
\end{figure}

In Figures \ref{fig:box} to \ref{fig:logistics}, red squares and gray dots indicate the locations of constructed facilities and affected areas respectively, small cyan nodes are the facility candidates that have not been established. The assignment between facilities and affected areas are connected by blue lines. The thickness of a line is proportion to the percentage of allocated demand, and the size of red squares represents the capacity of a facility.

We observe that all demands have been satisfied during the first three periods, and the demand allocation process gets finished at the end of day 2 under the penalty cost function $\mu_E(t)$. The reason is that $\mu_E(t)$ increases dramatically during the post-disaster phase, building more facilities with the largest available capacity can better serve the affected areas right after an emergency, and in turn decrease the total penalty cost. By contrast,  the penalty cost function $\mu_B(t)$ investigates a relatively small number of facilities during the initial hours and save the total logistic cost by enlarging facility capacities. The penalty cost function $\mu_L(t)$ establishes only 5 facilities out of 30 candidates, and dynamically adjust the corresponding capacities in periods 2 and 3. This is because $\mu_L(t)$ leads to a relatively larger unit penalty cost, satisfying demand during the initial hours will not significantly decrease the total penalty cost compared with the other two functions. Moreover, investigation of additional capacities is more economical than constructing new facilities. Thus, different types of unit penalty cost have different effect on the facility location and capacity planning decisions.
%, but the overall strategies are similar.


\subsection{Case study}\label{subsection:case}

\begin{figure}[htbp]
 \centering
 \includegraphics[width=0.6\linewidth]{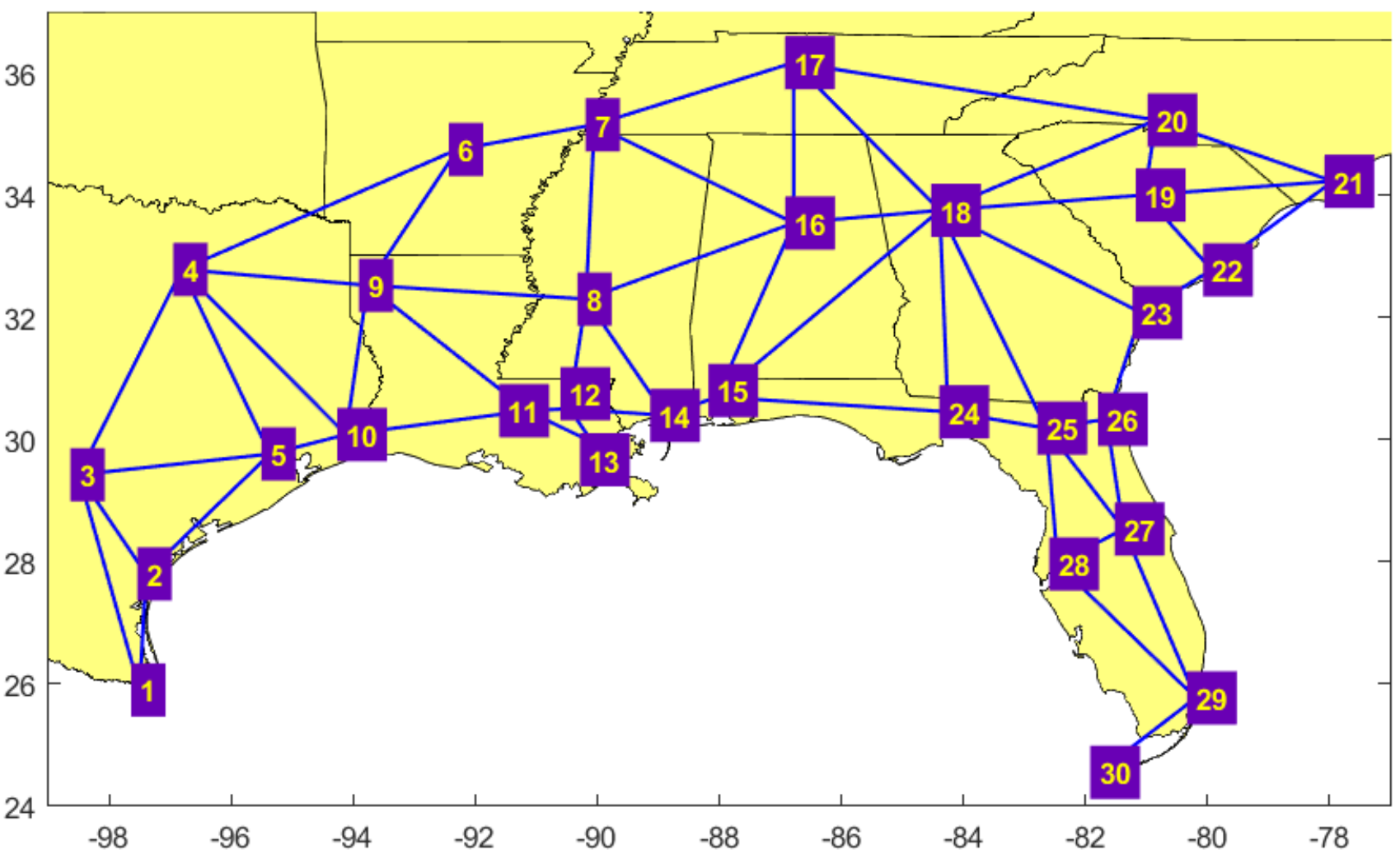}
 \caption{Case study network. The labeled number at each node is the index of each location.}
 \label{fig:case_label}
\end{figure}

The case study is based on a benchmark dataset proposed by \cite{rawls2010pre}, where a network with 30 nodes and 56 edges near the Gulf of Mexico under hurricane threats is incorporated, detailed data can be found in Figure \ref{fig:case_label}. Other input parameters, such as cost parameters and demand samples, are generally developed from \cite{rawls2010pre}, to make it more relevant, several modifications are purposed. Parameters of the case study are summarized in the following paragraph.

All nodes can be considered as affected areas and facility candidates, and $T$ is fixed as 3. During the preliminary hours right after a hurricane, it is difficult to establish large-scaled facilities, and the construction cost ($\bm f$) is inevitably expensive, thus, we assume that $\bm f=[300 000,188 400,94 200]$ and the maximum built capacity $\bm q=[204100,408200,780000]$ with respect to each time period. Varying construction cost of an additional capacity is randomly generated from 247.7 to 647.7, which is proportional to the purchase cost of new supplies. Unit penalty cost in each time period follows a DLF that is 500 times larger than that of Section \ref{subsection:performance}. As for the demand samples observed from historical data, we employ the same sample definition and occurrence probability as in Tables 3 and 4 of \cite{rawls2010pre}, and assume that the demand only arises from the landfall nodes when disaster occurs \citep{velasquez2020prepositioning}. For those hurricanes that do not have a specific landfall node, we uniformly divide the total demand to all the nodes. Details of the demand samples are summarized in Tables \ref{tab:scenario} and \ref{tab:demand} of \ref{app_scenario}.
%In a multi-period setting, we generally classify the stages of hurricanes into three stages. %,i.e., the preliminary stage right after the hurricane, the intermediate stage with steady and

\subsubsection{The location and capacity planning decisions in different periods}

\begin{figure}[htbp] \centering
\subfigure[$t=1$] { \label{fig:case1}
\includegraphics[width=0.45\linewidth]{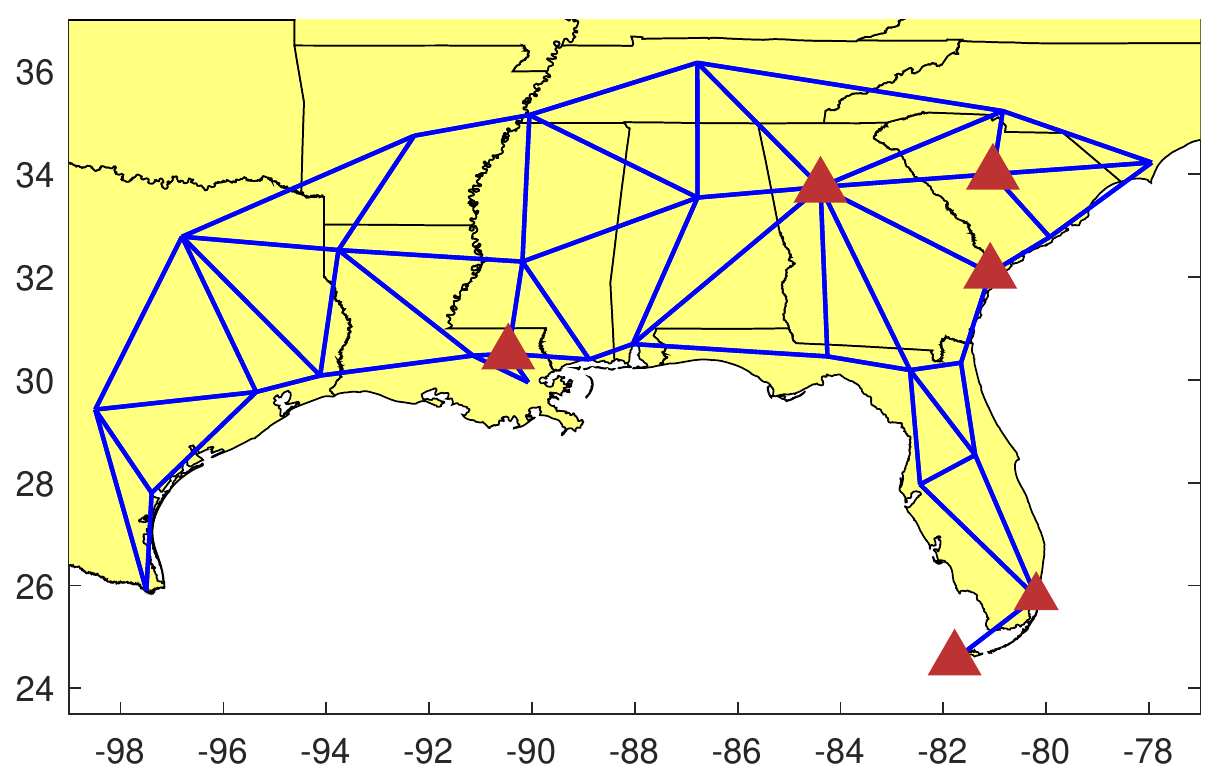}
}
\subfigure[$t=3$] { \label{fig:case2}
\includegraphics[width=0.45\linewidth]{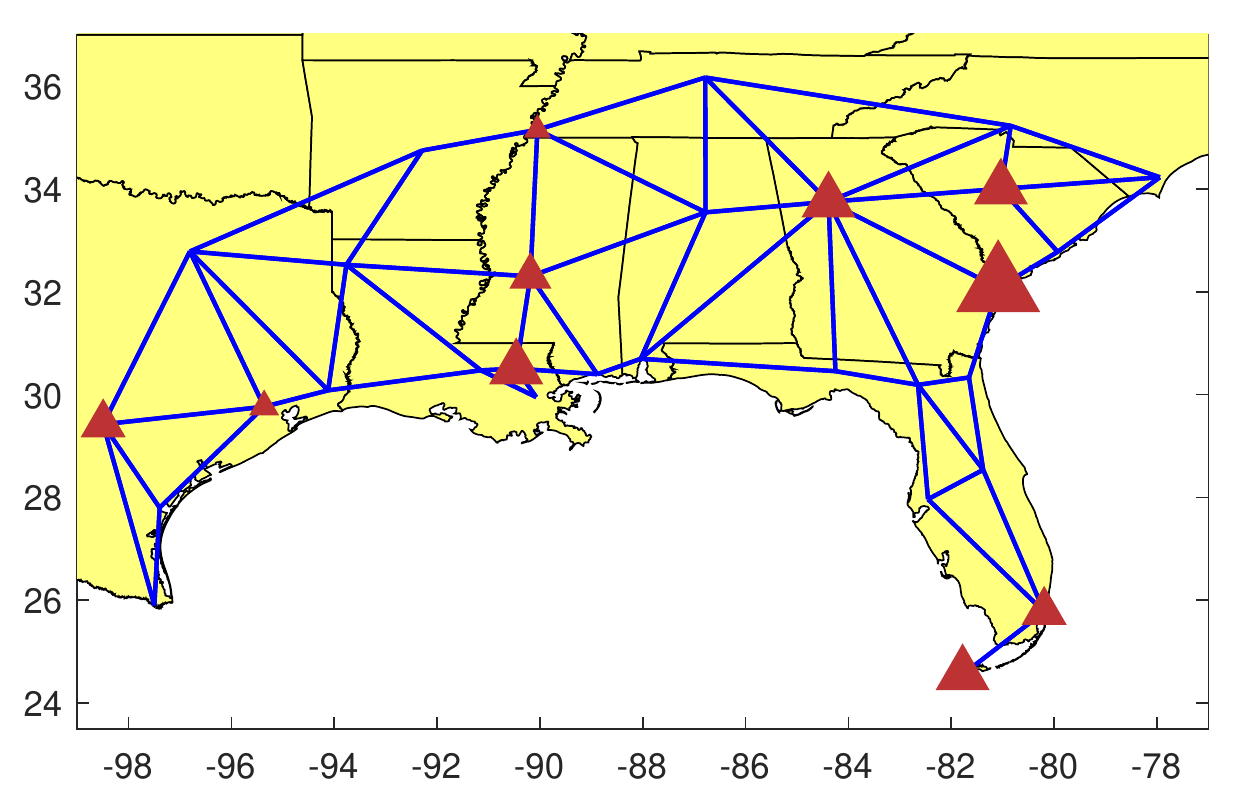}
}
\caption{The location and capacity planning decisions of the MFLCP-W model in different periods.}
\label{fig_case}
\end{figure}

Figure \ref{fig_case} illustrates the location and capacity planning decisions of the MFLCP-W model in different periods, where the size of red triangles is proportion to the volume of established capacity. During the preliminary stage of post-disaster phase, resources are scarce and construction costs are high, it is wise to establish a small amount of facilities with the maximum capacity at transportation hubs (such as node 18 Atlanta, GA) or non-disaster-prone locations closed to landfall nodes (such as node 12 Hammond, LA), which is validated by Figure \ref{fig:case1}. With time elapsing, there would be enough resource and personnel to build new facilities or expand the current facilities. Thus, more facilities are constructed, see the node 3 San Antonio, TX \& node 8 Jackson, MS in Figure \ref{fig:case2} and the capacity of some existing facilities are expanded, such as node 23 Savannah, GA in Figure \ref{fig:case2}.  Similar results can be found in the previous research using the same dataset, such as \cite{rawls2010pre} and \cite{velasquez2020prepositioning}.

%Thus, more facilities are constructed, see the node 3 San Antonio, TX \& node 8 Jackson, MS in Figure \ref{fig:case2} and the capacity of some existing facilities are expanded, such as node 23 Savannah, GA in Figure \ref{fig:case2}. 

%As depicted in Figure \ref{fig:case1}, during the preliminary stage of post-disaster phase, resources are scarce and construction costs are high. Therefore, it is wise to establish  facilities with the maximum capacity at transportation hubs (such as node 18 Atlanta, GA) or non-disaster-prone locations closed to landfall nodes (such as node 12 Hammond, LA). Similar results can be found in the previous research using the same dataset, such as \cite{rawls2010pre} and \cite{velasquez2020prepositioning}. With time elapsing, there would be enough resource and personnel to build new facilities at more locations or expand the current supply at some existing warehouses. Thus, more facilities are constructed, see the node 3 San Antonio, TX \& node 8 Jackson, MS in Figure \ref{fig:case2} and the capacity of some existing facilities are expanded, such as node 23 Savannah, GA in Figure \ref{fig:case2}. 

\subsubsection{Impact of the Wasserstein ball radius}

\begin{figure}[htbp] \centering
\subfigure[$H=30$] { \label{fig:case1_rad}
\includegraphics[width=0.45\linewidth]{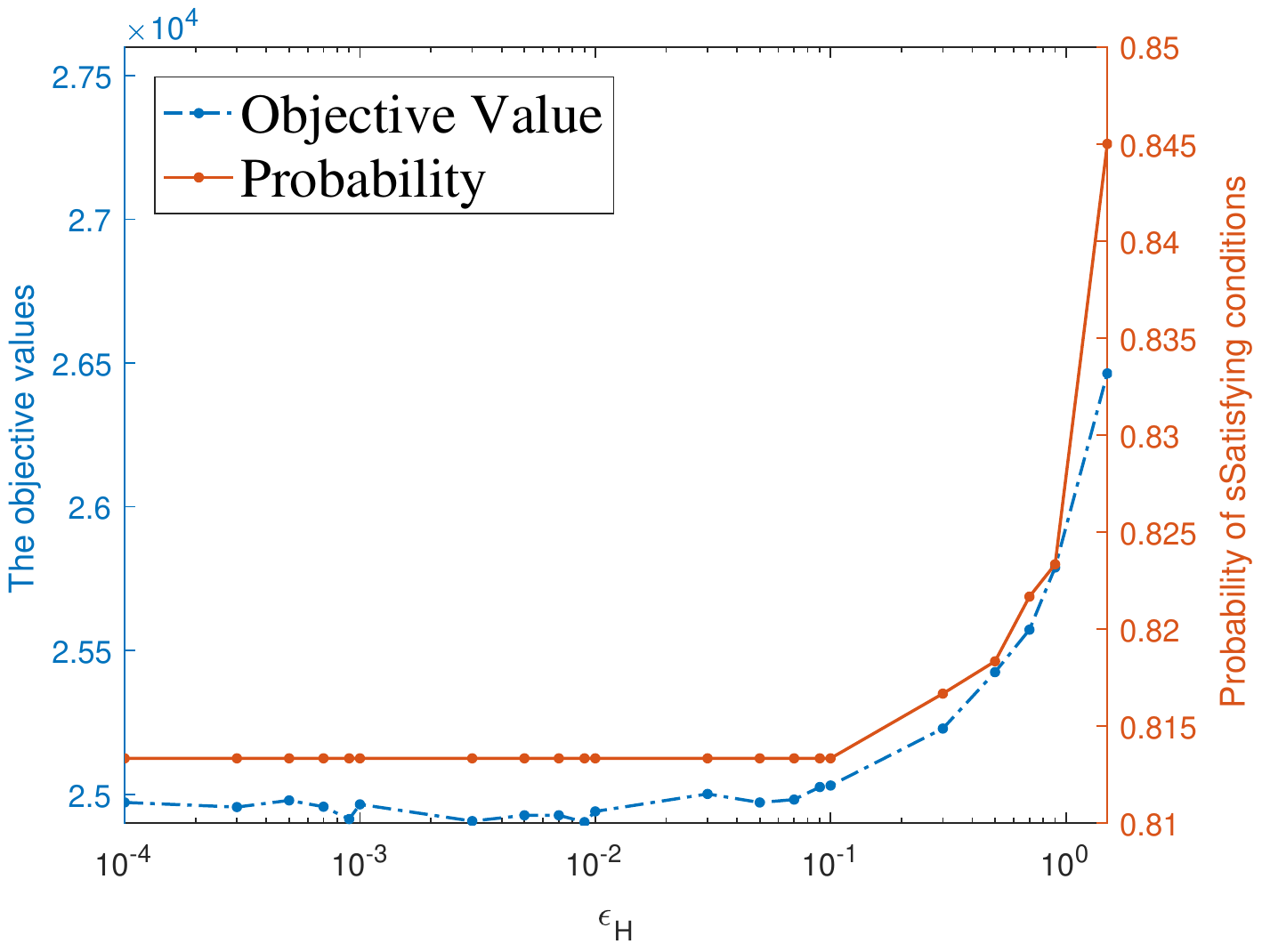}
}
\subfigure[$H=50$] { \label{fig:case2_rad}
\includegraphics[width=0.45\linewidth]{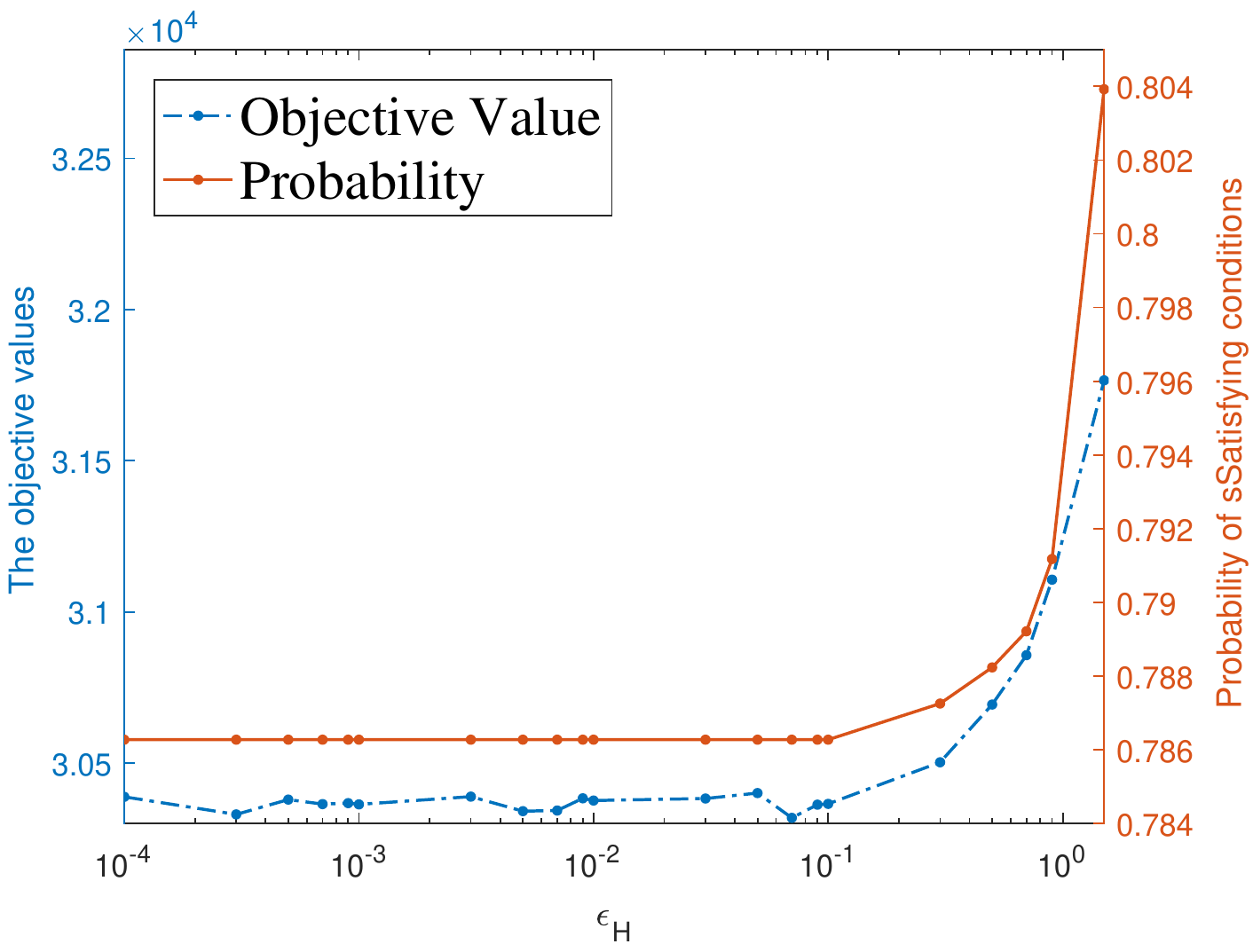}
}
\caption{The out-of-sample obj (left axis, blue line) and probability (right axis, orange line) as a function of the Wasserstein radius $\epsilon_H$ and averaged over $30$ experiments in the case study. (a) H = 30. (b)H = 50.}
\label{fig_case_rad}
\end{figure}

We conduct experiments to test the impact of radius $\epsilon_H$ on the MFLCP-W model over the benchmark dataset. The size of the training dataset is $H \in \{30,50\}$.  The averaged out-of-sample \ref{OBJ} and \ref{probability} over thirty independent experiments are shown in Figure \ref{fig_case_rad}. Similar to the results in Section \ref{subsection:radius}, the objective obtains its optimal value at a certain point, and then deteriorates with the increment of $\epsilon_H$; and the satisfaction probability is non-decreasing in $\epsilon_H$.

\subsubsection{Comparison with the state-of-the-art methods}

We conduct a performance comparison between the MFLCP-W model and the SAA model in this subsection.  As the maximum number of samples is 51 (see Table \ref{tab:scenario}), the sample size $H$ in this subsection varies from $10$ to $50$ with a step size of 5. We first use the same approach as Remark \ref{my_remark} to determine the optimal Wasserstein radius. Then, we compare the MFLCP-W model under the optimal Wasserstein radius with SAA in terms of the out-of-sample \ref{OBJ} and \ref{probability}. Results are summarized in Table \ref{tab:comparison},  where ``W" and ``SAA" represent the MFLCP-W and sample-based model respectively. Similar to the results in Section \ref{subsection:saa}, the proposed MFLCP-W model always achieves a higher satisfaction probability at a little sacrifice of operational cost. 

\begin{table}[htb]
 \centering
 \caption{Comparison with the SAA model}
 \begin{tabular}{cccccc}
 \toprule
  \multirow{2}[0]{*}{$H$} & \multirow{2}[0]{*}{Optimal Radius of W} & \multicolumn{2}{c}{Ratios of OBJ (W/S)} & \multicolumn{2}{c}{Probability} \\
  \cmidrule(r){3-4}\cmidrule(r){5-6}
     &  & W  & SAA & W  & SAA   \\ \midrule
    10 & 0.6 & 1.0176 & 1  & \textbf{0.650}  & 0.628 \\
    15 & 0.6 & 1.0147 & 1  & \textbf{0.750}  & 0.731 \\
    20 & 0.5 & 1.0136 & 1  & \textbf{0.776}  & 0.766\\
    25 & 0.5 & 1.0139 & 1  & \textbf{0.810}  & 0.801 \\
    30 & 0.3 & 1.0117 & 1  & \textbf{0.812}  & 0.809 \\
    35 & 0.3 & 1.0107 & 1  & \textbf{0.819}  & 0.811 \\
    40 & 0.3 & 1.0105 & 1 & \textbf{0.835}  & 0.834 \\
    45 & 0.3 & 1.0099 & 1 & \textbf{0.854}  & 0.849 \\
    50 & 0.25 & 1.0104 & 1  & \textbf{0.836}  & 0.828 \\
  \bottomrule
 \end{tabular}%
 \label{tab:comparison}%
\end{table}%

\subsubsection{Comparison with the single-period FLCP-W model}

In this subsection, we compare our MFLCP-W model with the single-period FLCP-W (SFLCP-W) model, i.e., $T=1$. The optimal Wasserstein radii for the MFLCP-W and SFLCP-W are determined by Remark \ref{my_remark} respectively. Then we compare the two models in terms of the out-of-sample \ref{OBJ} and \ref{probability}. Results are summarized in Table \ref{tab:comparison_s}, where ``W" and ``S" represent the proposed MFLCP-W and SFLCP-W models respectively.
\begin{table}[htb]
 \centering
 \caption{Comparison with the SFLCP-W model}
 \begin{tabular}{cccccc}
 \toprule
  \multirow{2}[0]{*}{$H$} & \multicolumn{2}{c}{Ratios of OBJ (W/S)} & \multicolumn{2}{c}{Probability} \\
  \cmidrule(r){2-3}\cmidrule(r){4-5}
      & W  & S & W  & S   \\ \midrule
    10  & 0.9158 & 1  & \textbf{0.663}  & 0.459 \\
    15  & 0.9272 & 1  & \textbf{0.750}  & 0.523 \\
    20  & 0.9418 & 1  & \textbf{0.745}  & 0.510\\
    25  & 0.9387 & 1  & \textbf{0.779}  & 0.542 \\
    30 & 0.9281 & 1  & \textbf{0.803}  & 0.556 \\
    35 & 0.9357 & 1  & \textbf{0.827}  & 0.606 \\
    40 & 0.9340 & 1 & \textbf{0.817}  & 0.562 \\
    45 & 0.9294 & 1 & \textbf{0.861}  & 0.595 \\
    50 & 0.9401 & 1  & \textbf{0.842}  & 0.584 \\
  \bottomrule
 \end{tabular}%
 \label{tab:comparison_s}%
\end{table}%

Table \ref{tab:comparison_s} indicates that the proposed MFLCP-W model always achieves a better out-of-sample performance than SFLCP-W in terms of the objective value and satisfaction probability, which confirms the advantages of MFLCP-W as expected.

\section{Conclusions}\label{section:conclusion}
In this paper, we propose a novel dat a-driven MFLCP-W model over the $\infty$-Wasserstein joint chance constraints to dynamically adjust facility locations and their capacities and guarantee a high probability for the on-time delivery. We reformulate the proposed MFLCP-W model as an MISOCP and design a tailored OA algorithm with provable convergence to solve it. Numerical results on the case study validate the better performance of our MFLCP-W model over its counterpart in terms of the cost and service quality. Moreover, the proposed OA algorithm significantly outperforms the commercial solver CPLEX 12.8 in terms of the computational time and optimality gap.   

Several future directions can be developed based on this research. First, other properties, such as equity and response time, can be incorporated to generate more managerial insights. Second, besides of opening new facilities, we can further decrease total logistics cost by closing and reopening facilities. Finally, it is also interesting to develop other algorithms to solve larger-scale problems.

%multi-period facility location and capacity planning problem during the PD-HL, where demand uncertainty is characterized by a $\infty$-Wasserstein ambiguity set. Joint chance constraints are incorporated to guarantee the probability of providing sufficiency reliefs to all victims is larger than an exogenous level. We employ a novel penalty cost function to better quantify victims' sufferings, which is monotone increasing with the deprivation time of relief (such as water, food or medicine).  After an equivalent transformation, the proposed DR model is reformulated as an MISOCP, and solved by a tailored OA algorithm with finite convergence. Numerical results revel a better out-of-sample performance of our model compared with the state-of-the-art SAA approach, find an powerful method to calibrate the Wasserstein radius, prove the asymptotic consistency of our data-driven approach, and validate the efficiency of the OA algorithm. We also investigate the sensitivity analysis of different types of demand-side cost. A practical case study under the hurricane threat is conducted in our experiments. The advantages the proposed model and several managerial insights are highlighted through these comprehensive numerical results.

\section*{Acknowledgement}
This research is partly supported by the National Natural Science Foundation of China (grant numbers 72101021, 72101012), the Fundamental Research Funds for the Central Universities 2021RC202, and the fellowship of China Postdoctoral Science Foundation (grant numbers 2021M690009 and 2021M690340).

\appendix
\section{Details of the hurricane samples}\label{app_scenario}

\cite{rawls2010pre} listed an instance of 15 historical hurricanes that has attacked the Atlantic Basin (East and Gulf Coast). The samples are classified into 5 categories, where categories 3 - 5 are characterized as major hurricanes, while the other two are minor ones. Potential demands and damage can be derived from the category of a hurricane. Table \ref{tab:demand} reports the demand vectors of each hurricane. After that, a set of 51 hurricanes samples, including both single and combined ones, are generated, and we estimate the corresponding probabilities according to the historical data. See Table \ref{tab:scenario} for more details.

% Table generated by Excel2LaTeX from sheet 'Sheet2'
\begin{table}[htbp]
 \centering\scriptsize
 \caption{Demand volume of each hurricane}
 \begin{tabular}{rcccccccccccccccc}
 \toprule
 \multicolumn{2}{c}{\diagbox{Node}{Hurricane} } & 1  & 2  & 3  & 4  & 5  & 6  & 7  & 8  & 9  & 10 & 11 & 12 & 13 & 14 & 15 \\\midrule
 1  & Brownsville & 0  & 0  & 0  & 0  & 0  & 0  & 0  & 0  & 0  & 75 & 0  & 600 & 94 & 0  & 0 \\
 2  & Corpus Christi & 0  & 0  & 0  & 0  & 0  & 0  & 0  & 0  & 0  & 75 & 0  & 600 & 94 & 0  & 0 \\
 3  & san antonio & 0  & 0  & 0  & 0  & 0  & 0  & 0  & 0  & 0  & 75 & 0  & 600 & 94 & 0  & 0 \\
 4  & dallas & 0  & 0  & 0  & 0  & 0  & 0  & 0  & 0  & 0  & 75 & 0  & 600 & 94 & 0  & 0 \\
 5  & Houston & 350 & 0  & 0  & 0  & 0  & 0  & 0  & 0  & 0  & 75 & 0  & 600 & 94 & 0  & 0 \\
 6  & Little Rock & 0  & 0  & 0  & 0  & 0  & 0  & 0  & 0  & 0  & 75 & 0  & 600 & 94 & 0  & 0 \\
 7  & Memphis & 0  & 0  & 0  & 0  & 0  & 0  & 0  & 0  & 0  & 75 & 0  & 600 & 94 & 0  & 0 \\
 8  & Jackson & 0  & 0  & 0  & 0  & 0  & 0  & 0  & 0  & 0  & 75 & 0  & 600 & 94 & 0  & 0 \\
 9  & Shreveport & 0  & 0  & 0  & 0  & 0  & 0  & 0  & 0  & 0  & 75 & 0  & 600 & 94 & 0  & 0 \\
 10 & Beaumont & 0  & 0  & 0  & 0  & 0  & 0  & 0  & 0  & 0  & 75 & 0  & 600 & 94 & 0  & 0 \\
 11 & Baton Rouge & 0  & 0  & 0  & 0  & 7500 & 0  & 0  & 1500 & 0  & 75 & 0  & 600 & 94 & 0  & 0 \\
 12 & Hammond & 0  & 0  & 0  & 0  & 0  & 0  & 0  & 0  & 0  & 75 & 0  & 600 & 94 & 0  & 0 \\
 13 & New Orleans & 0  & 0  & 0  & 0  & 0  & 0  & 0  & 0  & 1040 & 75 & 0  & 600 & 94 & 0  & 0 \\
 14 & Biloxi & 0  & 560 & 0  & 0  & 0  & 0  & 0  & 0  & 0  & 75 & 0  & 600 & 94 & 2239 & 0 \\
 15 & Mobile & 0  & 0  & 0  & 0  & 0  & 1000 & 0  & 0  & 0  & 75 & 0  & 600 & 94 & 0  & 0 \\
 16 & Birmingham & 0  & 0  & 0  & 0  & 0  & 0  & 0  & 0  & 0  & 75 & 0  & 600 & 94 & 0  & 0 \\
 17 & Nashville & 0  & 0  & 0  & 0  & 0  & 0  & 0  & 0  & 0  & 75 & 0  & 600 & 94 & 0  & 0 \\
 18 & Atlanta & 0  & 0  & 0  & 0  & 0  & 0  & 0  & 0  & 0  & 75 & 0  & 600 & 94 & 0  & 0 \\
 19 & Columbia & 0  & 0  & 0  & 0  & 0  & 0  & 0  & 0  & 0  & 75 & 0  & 600 & 94 & 0  & 0 \\
 20 & Charlotte & 0  & 0  & 0  & 0  & 0  & 0  & 0  & 0  & 0  & 75 & 0  & 600 & 94 & 0  & 0 \\
 21 & Wilmington & 0  & 0  & 0  & 0  & 0  & 0  & 600 & 0  & 0  & 75 & 5000 & 600 & 94 & 0  & 0 \\
 22 & Charleston & 0  & 0  & 861 & 9000 & 0  & 0  & 0  & 0  & 0  & 75 & 0  & 600 & 94 & 0  & 4400 \\
 23 & Savannah & 0  & 0  & 0  & 0  & 0  & 0  & 0  & 0  & 0  & 75 & 0  & 600 & 94 & 0  & 0 \\
 24 & Tallahassee & 0  & 0  & 0  & 0  & 0  & 0  & 0  & 0  & 0  & 75 & 0  & 600 & 94 & 0  & 0 \\
 25 & Lake City & 0  & 0  & 0  & 0  & 0  & 0  & 0  & 0  & 0  & 75 & 0  & 600 & 94 & 0  & 0 \\
 26 & Jacksonville & 0  & 0  & 0  & 0  & 0  & 0  & 0  & 0  & 0  & 75 & 0  & 600 & 94 & 0  & 0 \\
 27 & Orlando & 0  & 0  & 0  & 0  & 0  & 0  & 0  & 0  & 0  & 75 & 0  & 600 & 94 & 0  & 0 \\
 28 & Tampa & 0  & 0  & 0  & 0  & 0  & 0  & 0  & 0  & 0  & 75 & 0  & 600 & 94 & 0  & 0 \\
 29 & Miami & 0  & 0  & 0  & 0  & 7500 & 0  & 0  & 0  & 1040 & 75 & 0  & 600 & 94 & 0  & 0 \\
 30 & Key West & 0  & 0  & 0  & 0  & 0  & 0  & 0  & 0  & 0  & 75 & 0  & 600 & 94 & 2239 & 0 \\\bottomrule%\midrule
  %\multicolumn{2}{c}{ Category } & 3  & 5  &2  & 2  & 4  & 3  & 2  & 1  & 5  & 2 & 3  & 3 & 3 & 4  & 4 \\
 \end{tabular}%
 \label{tab:demand}%
\end{table}%

\begin{table}[htbp]
 \centering\scriptsize
 \caption{Scenario definitions and probability of occurrence \citep{rawls2010pre}. }
 \begin{tabular}{ccccccc}
 \toprule
 \multicolumn{3}{c}{Single hurricane threat} & \multicolumn{4}{c}{Combined hurricane threat} \\\cmidrule(r){1-3}\cmidrule(r){4-7}
 Sample & Hurricane & Probability & Sample & \multicolumn{2}{c}{Hurricanes} & Probability \\
 \midrule
 1  & 1  & 0.02308 & 16 & 1  & 2  & 0.0046 \\
 2  & 5  & 0.05 & 17 & 1  & 4  & 0.0057 \\
 3  & 10 & 0.16167 & 18 & 1  & 7  & 0.0057 \\
 4  & 3  & 0.05363 & 19 & 10 & 2  & 0.006 \\
 5  & 2  & 0.00925 & 20 & 10 & 13 & 0.0261 \\
 6  & 12 & 0.03083 & 21 & 10 & 9  & 0.0125 \\
 7  & 13 & 0.1338 & 22 & 10 & 8  & 0.0094 \\
 8  & 4  & 0.05363 & 23 & 2  & 5  & 0.005 \\
 9  & 11 & 0.02295 & 24 & 2  & 6  & 0.0047 \\
 10 & 14 & 0.02295 & 25 & 12 & 1  & 0.0052 \\
 11 & 15 & 0.02295 & 26 & 12 & 3  & 0.0061 \\
 12 & 7  & 0.05363 & 27 & 12 & 2  & 0.0047 \\
 13 & 9  & 0.05 & 28 & 12 & 4  & 0.0061 \\
 14 & 8  & 0.0308 & 29 & 12 & 14 & 0.0052 \\
 15 & 6  & 0.03083 & 30 & 13 & 2  & 0.0057 \\
   &  &  & 31 & 13 & 8  & 0.0086 \\
   &  &  & 32 & 4  & 2  & 0.005 \\
   &  &  & 33 & 11 & 5  & 0.0056 \\
   &  &  & 34 & 11 & 12 & 0.0052 \\\cline{1-3}
 \multicolumn{1}{|c}{}  &  & \multicolumn{1}{c|}{} & 35 & 11 & 13 & 0.0075 \\
 \multicolumn{2}{|c}{\multirow{2}[1]{*}{\shortstack{Total probability of\\ single hurricane}}} & \multicolumn{1}{c|}{\multirow{2}[1]{*}{0.75}} & 36 & 11 & 7  & 0.0057 \\
 \multicolumn{2}{|c}{} & \multicolumn{1}{c|}{} & 37 & 14 & 3  & 0.0057 \\
 \multicolumn{1}{|c}{} &  & \multicolumn{1}{c|}{} & 38 & 14 & 6  & 0.0052 \\
 \multicolumn{2}{|c}{\multirow{2}[1]{*}{\shortstack{Total probability of\\ combined hurricanes}}} & \multicolumn{1}{c|}{\multirow{2}[1]{*}{0.25}} & 39 & 15 & 5  & 0.0056 \\
 \multicolumn{2}{|c}{} & \multicolumn{1}{c|}{} & 40 & 15 & 7  & 0.0057 \\
\cline{1-3}   &  &  & 41 & 15 & 13 & 0.0075 \\
   &  &  & 42 & 15 & 14 & 0.005 \\
   &  &  & 43 & 9  & 1  & 0.0056 \\
   &  &  & 44 & 9  & 14 & 0.0056 \\
   &  &  & 45 & 8  & 5  & 0.006 \\
   &  &  & 46 & 8  & 3  & 0.0061 \\
   &  &  & 47 & 8  & 7  & 0.0061 \\
   &  &  & 48 & 6  & 5  & 0.006 \\
   &  &  & 49 & 6  & 3  & 0.0061 \\
   &  &  & 50 & 6  & 7  & 0.0061 \\
   &  &  & 51 & \multicolumn{2}{c}{No hurricane} & 0.0174 \\\bottomrule
 \end{tabular}%
 \label{tab:scenario}%
\end{table}%

	{\newpage

		\small
		\linespread{1.0} \selectfont
\bibliographystyle{elsarticle-harv}
\bibliography{sample}
}

\end{document}